\documentclass{article}
\usepackage[utf8]{inputenc}
\usepackage{complexity}
\usepackage{todonotes}
\usepackage{amsthm}
\usepackage{tikz}
\usepackage{xspace}
\usepackage[hidelinks]{hyperref}
\usepackage{amsmath,amsfonts}
\usepackage{stmaryrd} 

\def\rotateclockwise#1{
  \newdimen\xrw
  \pgfextractx{\xrw}{#1}
  \newdimen\yrw
  \pgfextracty{\yrw}{#1}
  \pgfpoint{\yrw}{-\xrw}
}

\def\rotatecounterclockwise#1{
  \newdimen\xrcw
  \pgfextractx{\xrcw}{#1}
  \newdimen\yrcw
  \pgfextracty{\yrcw}{#1}
  \pgfpoint{-\yrcw}{\xrcw}
}

\def\outsidespacerpgfclockwise#1#2#3{
  \pgfpointscale{#3}{
    \rotateclockwise{
      \pgfpointnormalised{
        \pgfpointdiff{#1}{#2}}}}
}

\def\outsidespacerpgfcounterclockwise#1#2#3{
  \pgfpointscale{#3}{
    \rotatecounterclockwise{
      \pgfpointnormalised{
        \pgfpointdiff{#1}{#2}}}}
}

\def\outsidepgfclockwise#1#2#3{
  \pgfpointadd{#2}{\outsidespacerpgfclockwise{#1}{#2}{#3}}
}

\def\outsidepgfcounterclockwise#1#2#3{
  \pgfpointadd{#2}{\outsidespacerpgfcounterclockwise{#1}{#2}{#3}}
}

\def\outside#1#2#3{
  ($ (#2) ! #3 ! -90 : (#1) $)
}

\def\cornerpgf#1#2#3#4{
  \pgfextra{
    \pgfmathanglebetweenpoints{#2}{\outsidepgfcounterclockwise{#1}{#2}{#4}}
    \let\anglea\pgfmathresult
    \let\startangle\pgfmathresult

    \pgfmathanglebetweenpoints{#2}{\outsidepgfclockwise{#3}{#2}{#4}}
    \pgfmathparse{\pgfmathresult - \anglea}
    \pgfmathroundto{\pgfmathresult}
    \let\arcangle\pgfmathresult
    \ifthenelse{180=\arcangle \or 180<\arcangle}{
      \pgfmathparse{-360 + \arcangle}}{
      \pgfmathparse{\arcangle}}
    \let\deltaangle\pgfmathresult

    \newdimen\x
    \pgfextractx{\x}{\outsidepgfcounterclockwise{#1}{#2}{#4}}
    \newdimen\y
    \pgfextracty{\y}{\outsidepgfcounterclockwise{#1}{#2}{#4}}
  }
  -- (\x,\y) arc [start angle=\startangle, delta angle=\deltaangle, radius=#4]
}

\def\corner#1#2#3#4{
  \cornerpgf{\pgfpointanchor{#1}{center}}{\pgfpointanchor{#2}{center}}{\pgfpointanchor{#3}{center}}{#4}
}

\def\hedgeiii#1#2#3#4{
  \outside{#1}{#2}{#4} \corner{#1}{#2}{#3}{#4} \corner{#2}{#3}{#1}{#4} \corner{#3}{#1}{#2}{#4} -- cycle
}

\def\hedgeiiii#1#2#3#4#5{
  \outside{#1}{#2}{#5} \corner{#1}{#2}{#3}{#5} \corner{#2}{#3}{#4}{#5} \corner{#3}{#4}{#1}{#5} \corner{#4}{#1}{#2}{#5} -- cycle
}

\def\hedgem#1#2#3#4{
  
  \outside{#1}{#2}{#4}
  \pgfextra{
    \def\hgnodea{#1}
    \def\hgnodeb{#2}
  }
  foreach \c in {#3} {
    \corner{\hgnodea}{\hgnodeb}{\c}{#4}
    \pgfextra{
      \global\let\hgnodea\hgnodeb
      \global\let\hgnodeb\c
    }
  }
  \corner{\hgnodea}{\hgnodeb}{#1}{#4}
  \corner{\hgnodeb}{#1}{#2}{#4}
  -- cycle
}

\def\hedgeii#1#2#3{
  \hedgem{#1}{#2}{}{#3}
}

\usepackage{enumitem}
\usetikzlibrary{calc, positioning, fit, shapes.misc}

\newcommand{\strat}{\mathcal{S}}

\newcommand{\PS}{\lang{Paired~SAT}\xspace}
\newcommand{\hyp}{H}
\newcommand{\WS}{E}
\newcommand{\som}{V}

\newcommand{\hxf}{\hyp = (\som, \WS)}

\newcommand{\CNF}{\lang{CNF}\xspace}
\newcommand{\QBF}{\lang{QBF}\xspace}
\newcommand{\etal}{{\em et al.}\xspace}
\newcommand{\hedge}{edge\xspace}
\newcommand{\hedges}{edges\xspace}

\DeclareMathOperator{\CL}{CL}
\DeclareMathOperator{\CCL}{CCL}

\newcommand{\Client}{\mathcal{C}}
\newcommand{\Waiter}{\mathcal{W}}
\newcommand{\uCwin}{u\text{-}\mathcal{F}}
\newcommand{\vCwin}{v\text{-}\mathcal{F}}

\newtheorem{theorem}{Theorem}
\newtheorem{lemma}[theorem]{Lemma}
\newtheorem{corollary}[theorem]{Corollary}
\newtheorem{definition}[theorem]{Definition}
\newtheorem{fact}[theorem]{Fact}
\newtheorem{claim}[theorem]{Claim}

\newtheorem{proposition}[theorem]{Proposition}

\newcommand{\qedclaim}{\hfill $\diamond$ \medskip}

\newenvironment{proofclaimitem}{\noindent{\em Proof of the claim.}}{}

\tikzstyle{v}=[circle,inner sep=0, minimum size =6 pt, line width = 1pt, draw=black, fill=black, text= white]
\tikzstyle{he}=[rounded corners, draw=black]
\tikzstyle{noeud}=[circle, inner sep=2, minimum size =3 pt, line width = 1pt, draw=black, fill=white]

\title{On the complexity of Client-Waiter and Waiter-Client games
\thanks{This research was partly supported by the ANR project P-GASE (ANR-21-CE48-0001-01)}
}
\author{
Valentin Gledel, Nacim Oijid, Sébastien Tavenas, Stéphan Thomassé
}
\date{}

\begin{document}

\maketitle

\begin{abstract}
    Positional games were introduced by Hales and Jewett in 1963, and their study became more popular after Erd\H os and Selfridge's first result on their connection to Ramsey theory and hypergraph coloring in 1973. Several conventions of these games exist, and the most popular one, Maker-Breaker was proved to be \PSPACE-complete by Schaefer in 1978. The study of their complexity then stopped for decades, until 2017 when Bonnet, Jamain, and Saffidine proved that Maker-Breaker is $W[1]$-complete when parameterized by the number of moves. The study was then intensified when Rahman and Watson improved Schaefer's result in 2021 by proving that the \PSPACE-hardness holds for $6$-uniform hypergraphs. More recently, Galliot, Gravier, and Sivignon proved that computing the winner on rank~$3$ hypergraphs is in \P. 
    
    We focus here on the Client-Waiter and the Waiter-Client conventions. Both were proved to be \NP-hard by Csernenszky, Martin, and Pluhár in 2011, but neither completeness nor positive results were known for these conventions. In this paper, we complete the study of these conventions by proving that the former is \PSPACE-complete, even restricted to $6$-uniform hypergraphs, and by providing an \FPT-algorithm for the latter, parameterized by the size of its largest \hedge. In particular, the winner of Waiter-Client can be computed in polynomial time in $k$-uniform hypergraphs for any fixed integer $k$. Finally, in search of finding the exact bound between the polynomial result and the hardness result, we focused on the complexity of rank~$3$ hypergraphs in the Client-Waiter convention. We provide an algorithm that runs in polynomial time with an oracle in \NP.
\end{abstract}

\section{Introduction}

Positional games were introduced by Hales and Jewett~\cite{HJ63} as a generalization of the Tic-Tac-Toe game. These games are played on a hypergraph $H = (V,E)$, on which the players alternately claim an unclaimed vertex of $V$. 
The Tic-Tac-Toe corresponds to the convention Maker-Maker: the first player who claims all vertices of an \hedge of \(E\) wins. However, since the second player cannot win in Maker-Maker games, and since Maker-Maker games are not hereditary, for the outcome, the study quickly switched to another convention: Maker-Breaker.

In the most studied convention, Maker-Breaker, Maker (one player) tries to claim all the vertices of an \hedge, while Breaker (the other player) aims to prevent her from doing so. If not explicitly specified, we assume here that Maker plays the first move. The study became more popular in 1973 when Erd\H os and Selfridge~\cite{ES73} obtained the following criterion.
\begin{theorem}[Erd\H os-Selfridge criterion]
\label{thm:ErdosSelfridge}
Let \(H = (V,E)\) be a hypergraph. If
\[
    \sum_{e \in E} 2^{-\lvert e\rvert} < \frac{1}{2},
\]
then the position is a win for Breaker. 
\end{theorem}

For the biased versions of Maker-Breaker, for example in the \((1 : b)\) one (at each turn Maker selects one vertex then Breaker selects \(b\) of them), one can ask for different families of hypergraphs, what is the threshold for \(b\) which turn the position for Breaker from loser to winner. Surprisingly, for a large number of hypergraphs families (but not all), it has been found that the threshold for the random game (where both players play randomly) is essentially the same as the threshold for the optimal play. This phenomenon, known as the "probabilistic intuition" was first demonstrated by Chvatal and Erd\H os~\cite{CE78}, and then further developed in numerous papers particularly by Beck (see for example the books~\cite{Bec02} and~\cite{HKST14}).



Schaefer proved in 1978 that it is \PSPACE-complete to determine the winner of a Maker-Breaker game~\cite{Sch78} (the problem appears there under the name \(G_{\text{pos}}(\text{POS CNF})\)). More precisely, he shows the \PSPACE-hardness even for hypergraphs of rank at most \(11\) (i.e., such that the size of the \hedges is bounded by \(11\)). A simplified proof can be found in~\cite{Bys04}. The result was improved in 2021 by Rahman and Watson~\cite{RW21}: the problem is \PSPACE-complete for $k$-uniform hypergraphs with \(k \ge 6\) (i.e., hypergraphs where all \hedges have size exactly \(k\)). 
Byskoz~\cite{Bys04} also notices that deciding the winner in a Maker-Breaker can be reduced to the same problem in the Maker-Maker convention (up-to increasing by \(1\) the maximal size of its edges). In particular, deciding who wins in a Maker-Maker game is \PSPACE-complete for \(k\)-uniform hypergraphs as soon as \(k \ge 7\).
 
On the positive side, Kutz~\cite{Kut05} proved in 2005 that the problem is tractable for $2$-uniform hypergraphs and $3$-uniform linear ones\footnote{Hypergraphs whose intersection of each pair of \hedges is of size at most one.}. It was improved recently by Galliot \etal~\cite{Gal23thesis, GGS22}: deciding the winner is tractable for rank~$3$ hypergraphs.

Since the introduction of positional games, the studies have focused on the Maker-Breaker convention. In order to have a better understanding of this convention, Beck~\cite{Bec02} introduced the Client-Waiter and Waiter-Client conventions in 2002 under the names Chooser-Picker and Picker-Chooser. Their current names were suggested by Hefetz, Krivelevich, and Tan~\cite{HKT16} as they are less ambiguous. In both conventions, Waiter selects two vertices of the hypergraph and offers them to Client. Client then chooses one to claim, and the second one is given to Waiter. If the number of vertices is odd, the last vertex goes to Client. In the Client-Waiter convention, Client wins if he claims all vertices of an \hedge, otherwise Waiter wins. In the Waiter-Client convention, Waiter wins if she claims all vertices of an \hedge, otherwise Client wins. 

Notice that these conventions can also be seen as variations of  Avoider-Enforcer. In particular, the definition we gives for Waiter-Client is the one which appears originally in Beck~\cite{Bec02}. But since~\cite{HKT16}, Waiter-Client is often defined following the Avoider-Enforcer convention: Waiter wins if she can forces Client to claim a whole edge (and if the number of vertices is odd, the last vertex goes to Waiter). One can notice that both definitions correspond in fact to the same game. Unlike the fact that Maker-Breaker and Avoider-Enforcer are very different games, in this ``I-cut-you’ll-choose way'' paradigm, since Client picks a vertex from only two possibilities, it is symmetric to associate the chosen vertex to Client or to Waiter.

The study of Client-Waiter and Waiter-Client games were first motivated by its similarities to Maker-Breaker games. 
For example, Bednarska-Bzd\c{e}ga (improving previous results~\cite{Bec02,CMP09}) proved that Theorem~\ref{thm:ErdosSelfridge} also holds in Waiter-Client convention~\cite{Bed13}. Moreover, the ``probabilistic intuition'' continues to work~\cite{Bec02} in these conventions. More recently, this probabilistic method has been stated in a more general case for the biased Waiter-Client $H$-game by Bednarska-Bzd\c ega, Hefetz and \L{}uczak in 2016~\cite{Bed16}. This conjecture has just been proved recently by Nenadov in 2023~\cite{Nen23}. These similarities led Beck~\cite{Bec02} and Csernenszky, Mándity, and Pluhár~\cite{CMP09} to conjecture that if Maker wins in a Maker-Breaker game on some hypergraph $H$ going second, then Waiter wins in Waiter-Client. Notice that up to considering the transversal of \(H\), the conjecture also implies that a win for Breaker as a second player implies a win for Waiter in the Client-Waiter game. But this conjecture was disproved by Knox~\cite{Kno12} in 2012. Today however, Client-Waiter and Waiter-Client games are studied independently of Maker-Breaker games: Csernenszky~\cite{Cse10} proved in 2010 that $7$-in-a-row Waiter-Client is a Client win on the infinite grid, while this problem is still open in the Maker-Breaker convention, and Hefetz \etal~\cite{HKT16} studied several classical games in Waiter-Client convention.

In terms of complexity, only few results were known about Waiter-Client and Client-Waiter games. In contrast to Maker-Breaker, Maker-Maker, Avoider-Avoider and Avoider-Enforcer, which are known to be \PSPACE-complete~\cite{Sch78, Bys04, BH19, RW21, GO23}, the asymmetry between the players moves in Waiter-Client and Client-Waiter conventions makes it more difficult to obtain reductions. In fact, Waiter has more choices than Client in her moves, which makes most reduction techniques fail. Both problems have been conjectured \PSPACE-complete in~\cite{CMP09}. In~\cite{CMP11}, the authors show that both problems are \NP-hard. However, in the Waiter-Client convention, the hypergraph has an exponential number of \hedges but is given in a succint way: it is the transversal of a given hypergraph. This leaves open the question of whether the problem of deciding who wins a Waiter-Client game is still \NP-hard when the hypergraph is given via the list of its edges.

The study of the parameterized complexity of combinatorial games emerged roughly together with the study of parameterized problems, and some strong results about complexity theory are due to this study. For instance, Abrahamson, Rodney, Downey, and Fellows\cite{Abr93} proved that $AW[3] = AW[*]$, through the game Geography.
Only few results are known
on positional games with the parameterized complexity paradigm, and only related to three conventions: Maker-Breaker, Maker-Maker and Avoider-Enforcer. Their study started by Downey and Fellows, conjecturing that {\sc Short Generalized Hex} was \FPT, but it was disproved (unless \FPT = \W[1]) by Bonnet, Jamain and Saffidine in 2016~\cite{Bon16}, proving its \W[1]-hardness. Then, in 2017, Bonnet \etal~\cite{Bon17} proved that general Maker-Breaker games are \W[1]-complete, Avoider-Enforcer games are co-\W[1]-complete and general Maker-Maker games are \AW[*]-complete.

Notice that the parameterized results obtained on general positional game consider the number of moves as a parameter. This is mostly motivated by the fact that these problems are already \PSPACE-hard for bounded rank hypergraphs. In Waiter-Client convention however, this is not the case, and therefore, we study the complexity of determining the winner of Waiter-Client games parameterized by the rank of the hypergraph. Note that, even if the outcome of Client-Waiter games are hereditary (see Lemma~\ref{lem: sub-hypergraph}), in contrast with Maker-Breaker or Avoider-Enforcer games, when the number of moves is bounded, it is not. Indeed, Waiter can control where Client plays, the addition of isolated vertices can be used by Waiter to waste turns. Therefore, the number of moves in this convention is not as relevant as in the others.

\subsubsection*{High-level description of the results.} 

We show that, similarly to the Maker-Breaker and Enforcer-Avoider conventions, deciding the winner of a positional game in the Client-Waiter convention is \PSPACE-complete. The result was already conjectured in 2009~\cite{CMP09}. It is an improvement of~\cite{CMP11} where the problem is shown to be \NP-hard. Moreover, we obtain the \PSPACE-completeness  even for \(6\)-uniform hypergraphs.

\begin{theorem}\label{6-uniform Client-Waiter}
    For \(k \ge 6\), Client-Waiter games are \PSPACE-complete even restricted to $k$-uniform hypergraphs.
\end{theorem}

The containment in \PSPACE~directly follows from Lemma~2.2 in~\cite{Sch78}. So the main point is the hardness part. To prove it, we reduce the problem of deciding who has the win to the same problem in the following game. 

\begin{definition}[\PS] \label{def: paired-sat}
    Let $\varphi$ be a $3$-CNF Formula over a set of pairs of variables $X = \{(x_1, y_1), \dots, (x_n, y_n)\}$. The \PS-game is played by two players, Satisfier and Falsifier as follows: while there is a variable that has not been assigned a valuation, Satisfier chooses a pair of variables $(x_i,y_i)$ that she has not chosen yet and gives a valuation, $\top$ or $\bot$, to $x_i$. Then Falsifier gives a valuation to $y_i$. When all variables are instantiated, Satisfier wins if and only if the valuation they have provided to the $x_i$s and $y_i$s satisfies $\varphi$.
\end{definition}

The \PS-game is a new variant of a \CNF-game where the play order is closer to the Client-Waiter convention: the first player chooses, at each turn, which variables she plays on and which variable the second player will have to play on.

Again, deciding who wins on this new game is \PSPACE-complete. It is proved in Section~\ref{sec:PairedSAT} by a reduction from the game \(3\)-QBF (known to be \PSPACE-complete since Schaefer's seminal work~\cite{Sch78}).

\begin{theorem}\label{thm:PS-PSPACE}
    Deciding who is the winner of the \PS-game is \PSPACE-complete.
\end{theorem}

Then Theorem~\ref{6-uniform Client-Waiter} is obtained by reducing the \PS-game to the Client-Waiter one. This is done by designing a gadget (given in Figure~\ref{fig: gadget client waiter}) which simulates a pair of variables \((x_i,y_i)\) of the \PS-game. Notice that the reduction only creates a hypergraph of rank \(6\). But, similarly to the Maker-Breaker games~\cite[Corollary~4]{RW21} and the Avoider-Enforcer ones~\cite[Lemma~7]{GO23}, Lemma~\ref{to 6uniform Client-Waiter} shows that the hypergraph can be turned into a \(6\)-uniform one afterwards. 

\medskip

In the Maker-Breaker convention, as said before, it is known~\cite{RW21} that deciding if a position is winning is \PSPACE-complete over hypergraphs of rank at most \(6\). On the other side, the problem is easy for hypergraphs of rank 2 since Maker wins if and only if there are two adjacent \(2\)-edges or the graph contains a singleton \hedge (result already noticed in~\cite{Kut05}). But only recently, after a serie of results~\cite{Kut05, RW20}, Galliot, Gravier, and Savignon~\cite{Gal23thesis, GGS22} showed that the problem is still polynomial for hypergraphs of rank at most \(3\).
The same question arises in the Client-Waiter convention: despite being \PSPACE-complete over hypergraphs of rank at most \(6\), what can we say about the complexity of the problem for hypergraphs of low rank? For hypergraphs of rank \(2\), it is easily seen that the problem is polynomial (Proposition~\ref{prop: rank 2 Client-Waiter}). The question is already non-trivial for hypergraphs of rank \(3\). 

We show that the problem of deciding if a position is winning in a rank~\(3\) hypergraph reduces to the problem of finding a specific structure (called a Tadpole) in a hypergraph.
Tadpoles have been generalized to hypergraphs by Galliot, Gravier and Sivignon~\cite{GGS22} to handle rank~3 Maker-Breaker games and their definition is recalled in Definition~\ref{def_tadpole}. Intuitively, given \(a\) and \(b\) two vertices of \(H\), an \(ab\)-{\em tadpole} is given by the union of a path from \(a\) to \(b\) and a cycle containing \(b\) such that the structure is simple (two \hedges intersect only if they are two consecutive \hedges of the path or the cycle, or if one is the last \hedge of the path and the other an \hedge of the cycle containing \(b\)) and linear (the intersection of two intersecting \hedges is of size exactly~\(1\)). More precisely, we require this structure is \(3\)-uniform, simple and linear, which means that all edges have size \(3\), the intersection of two consecutive edges has always size one, and a same vertex can not happen at two different places of the structure. An \(ab\) tadpole is simply called {\em tadpole}. A tadpole is said rooted in \(a\), if it is an \(ab\)-tadpole for some vertex \(b\).

We consider the problem \(\textsc{Tadpole}\):  Given a vertex \(a\) in a \(3\)-uniform hypergraph \(H\), decide if there is in \(H\) a tadpole rooted in \(a\).

It is easy to check that a given subhypergraph is a tadpole, so this problem is in \(\NP\). We do not know if this problem can be tractable. We note however that it would be sufficient to loop for all vertex \(b\) and all triples of edges of the form \(\{b,x_1,x_2\},\{b,y_1,y_2\},\{b,z_1,z_2\}\) and check if there are two {\em disjoint} simple linear paths linking the two sources \(a\) and \(x_1\) to the targets \(y_1\) and \(z_1\) (\(b\) would be the contact between the path and the cycle). The complexity of the problem ``Disjoint Connected Paths'' for graphs has been a very fruitful research topic. In the case of two sources and two targets, the problem was shown to be tractable~\cite{Shi80,Sey80}. In fact, in their well-known result, Robertson and Seymour~\cite{RS95} showed that the problem continues to be tractable for a constant number of sources and targets. The case where the number of sources and targets is unbounded is one of the first \NP-complete problems in Karp's list. For \(3\)-uniform hypergraphs, it has been just proved recently~\cite{GGS23} that the simple linear connectivity problem (one source and one target) is tractable. A corollary of the next theorem is that if the ``Disjoint Connected Paths'' problem for two sources and two targets is still tractable for \(3\)-uniform hypergraphs, then deciding who is winning in a Client-Waiter game on a hypergraph of rank \(3\) would also be tractable.

\begin{theorem}\label{thm: rank3 client-waiter}
	The problem of deciding if a given rank~\(3\) hypergraph is a winning position for Client can be solved by a polynomial time algorithm which uses the problem \(\textsc{Tadpole}\) as an oracle.
	
	In particular, the problem lies in the class $\Delta_2^{\P} = \P^{\NP}$.
\end{theorem}

In fact, during the proof of this theorem, we notice that this reduction is necessary. The problem of detecting a Tadpole can conversely be reduced to deciding on a winning position in a Client-Waiter game.

\begin{proposition}\label{prop:HarderThanTadpole}
	The problem \textsc{Tadpole} is polynomial-time many-one reducible to the problem of deciding if Client has a winning strategy in a Client-Waiter game played on a rank~$3$ hypergraph.
\end{proposition}

\medskip

Then, we focus on the second convention Waiter-Client. Surprisingly, the complexity of deciding if a position is winning is very different in this convention. 

\begin{theorem}\label{th:FPT}
    Waiter-Client is \FPT\ on $k$-uniform hypergraphs. 
    
    More precisely deciding if a hypergraph \(H=(V,E)\) is winning for Waiter can be decided in time \(O(f(k) \lvert E \rvert \log \lvert V \rvert)\) where \(f\) is a computable function, i.e., in linear time when \(k\) is fixed.
\end{theorem}

This result is obtained from a structural analysis of $k$-uniform hypergraphs, using the famous sunflower lemma from Erd\H os and Rado~\cite{ER60}. This is a very natural approach since, intuitively, a huge sunflower (\hedges pairwise intersecting on a same center set) should be ``reducible'' in the sense that one could replace the sunflower by a unique \hedge given by the center set. Indeed, if Waiter can obtain the center and the sunflower is large enough, then she can make sure to get a petal and so a set of the sunflower. This intuition should be valid, but the main difficulty is to exactly state what is ``huge'', as all other potentially useful sunflowers interact. We were unable to find a simple strategy based on these lines.
Instead, we describe a kernelization type algorithm which produces a subset of vertices (kernel) of the hypergraph $H$ such that Waiter wins on $H$ if and only if Waiter wins on the trace of $H$ on the kernel. The argument is based on a process which ultimately reaches a fixed point, the major drawback being that the kernel size is ridiculously large. 
This kernelization provides an \FPT\ algorithm for the Waiter-Client game on rank~$k$ hypergraphs. This also proves that if Waiter can win, she wins using only some function of $k$ moves. However, the gap between the very large upper bound and the best known lower bound ($2^k-1$ moves, required to win on $2^k$ disjoint \hedges of size $k$) indicates how little we understand strategies in Waiter-Client games.

Nevertheless, the fact that Waiter-Client is FPT for rank~$k$ hypergraphs could indicate that this convention is simpler to analyse than Maker-Breaker. It would be interesting to revisit the classical topics in which Maker-Breaker game was tried as a tool (for instance 2-colorability of hypergraphs or the Local Lemma) to see if Waiter-Client could provide more insight.


Complexity results for the various conventions are summarized in the table below.

\begin{table}[ht]
    \centering
\begin{tabular}{  l || p{1.55cm} | p{1.55cm} | p{1.55cm}  |p{1.55cm}  | p{2.25cm}}
    Rank r & 2 & 3 & 4, 5 & 6 & 7+ \\
    \hline  \hline 
   Maker &  P[Folklore]   &  P\cite{Gal23thesis}   & Open & \PSPACE-c & \PSPACE-c \\
   ~-Breaker &  &  & & \cite{RW21}  & \cite{RW21, Sch78}  \\ & & & & & \\ \hline 
   Maker & P[Folklore]  & Open & Open & Open & \PSPACE-c \\
   ~-Maker & & & & & \cite{RW21, Bys04} \\ & & & & & \\  \hline 
   Avoider & \PSPACE-c & \PSPACE-c & \PSPACE-c & \PSPACE-c &  \PSPACE-c \\ 
   ~-Avoider & \cite{BH19} & \cite{BH19}  & \cite{BH19} & \cite{BH19} & \cite{BH19} \\ & & & & & \\ \hline 
   Avoider & P~\cite{Gal23thesis} & Open & Open & \PSPACE-c & \PSPACE-c \\ 
   ~-Enforcer & & & & \cite{GO23}  & \cite{GO23} \\ & & & & & \\ \hline 
   Client &  P Prop~\ref{prop: rank 2 Client-Waiter} & $\Delta_2^{\P} = \P^{\NP}$ & Open & \PSPACE-c & \PSPACE-c \\
   ~-Waiter & & Cor~\ref{thm: rank3 client-waiter} & & Thm~\ref{6-uniform Client-Waiter} & Thm~\ref{6-uniform Client-Waiter} \\ & & & & & \\ \hline 
   Waiter & P Prop~\ref{prop: rank 2 Waiter-Client} & P Thm~\ref{th:FPT} & P Thm~\ref{th:FPT} & P Thm~\ref{th:FPT} & FPT w.r.t.\ $r$ \\
   ~-Client & & & & & Thm~\ref{th:FPT}
 \end{tabular}
    \caption{Complexity in the different conventions}
    \label{tab: complexity of conventions}
\end{table}

\subsubsection*{Organization.} 

We start with some preliminary results and definition in Section~\ref{sec: preliminaries}. We then prove in Section~\ref{sec: 6unif client-waiter} that Client-Waiter games are \PSPACE-complete, even restricted to $6$-uniform hypergraphs. In Section~\ref{sec: rank3 client-waiter}, we focus on rank~3 hypergraph in Client-Waiter convention and we prove that the decision problem of determining the outcome of the game is in $\Delta_2^{\P}$. Finally, in Section~\ref{sec: waiter-client}, we prove that Waiter-Client games are \FPT\ parameterized by the rank.

\section{Preliminaries}\label{sec: preliminaries}

In this section, we first introduce the context of the game by providing some definitions, then we present some useful lemmas to handle Client-Waiter and Waiter-Client games.

\begin{definition}
    Let $H = (V,E)$ be a hypergraph and let $k$ be an integer. $H$ is said to have {\em rank~$k$} if all its \hedges $e \in E$ have size at most $k$. $H$ is said to be $k$-uniform if all its \hedges have size exactly $k$.
\end{definition}

Remark that, if $H$ is a hypergraph, if it has an \hedge included in another, we can remove the largest one without changing the outcome of the game played on $H$. Therefore, we can consider that all hypergraphs considered in this paper are clutter.

\begin{definition}
    Let $H = (V,E)$ be a hypergraph. $H$ is said to be a {\em clutter} if for any \hedges $e_1 \neq e_2 \in E$, we have $e_1 \not \subseteq e_2$ and $e_2 \not \subseteq e_1$.
\end{definition}

\begin{definition}
    The hypergraph \(H'=(V',E')\) is a subhypergraph of \(H = (V,E)\) if \(V' \subseteq V\) and \(E' \subseteq E\). If \(A \subseteq V\), the induced subhypergraph \(H_{|A}\) is the subhypergraph \((A,\{e \in E \mid e \subseteq A\})\). The trace \(T_A(H)\) of \(H\) on \(A\) is the hypergraph \((A,\{e \cap A \mid e \in E\})\). 
\end{definition}

In this section, we present some general results about Client-Waiter and Waiter-Client games. 
The monotony of Client-Waiter and Waiter-Client conventions is immediate and folkloric: 
if ``Maker'' player (i.e.\ Client in Client-Waiter and Waiter in Waiter-Client) has a winning strategy on a sub-hypergraph, then he has one on the general hypergraph.
\begin{lemma}\label{lem: sub-hypergraph}
    Let $\hxf$ be a hypergraph and let $\hyp' = (\som', \WS')$ be a sub-hypergraph of $\hyp$. If Client (resp. Waiter) wins the Client-Waiter (resp. Waiter-Client) game on $\hyp'$, then Client (resp. Waiter) wins the Client-Waiter (resp. Waiter-Client) game on $\hyp$.
\end{lemma}

In the Client-Waiter convention, \hedges of length \(2\) are forced moves for Waiter.

\begin{lemma}[Proposition 9 from \cite{CMP09}]
\label{lem:size2hyperedge}
Let $\hxf$ be a hypergraph. Let a game in the Client-Waiter convention. Let $W,C \subset \som$ be the set of vertices already claimed by Waiter and Client respectively. If there exists $e\in \WS$ such that $e \cap W = \emptyset$ and $|e \setminus C| = 2$, then an optimal move for Waiter is to propose the two unclaimed element of $e$ with her next move.
\end{lemma}

\section{6-uniform Client-Waiter games are \PSPACE-complete} \label{sec: 6unif client-waiter}

This section is dedicated to the proof of Theorem~\ref{6-uniform Client-Waiter}. As explained in the introduction, we start with hypergraphs of rank at most \(6\):
\begin{proposition}\label{prop: Client-Waiter pspace complete}
    Computing the winner of a Client-Waiter game is \PSPACE-complete, even restricted to hypergraphs of rank~$6$.
\end{proposition}

We notice that the membership in \PSPACE\ follows from an argument of Schaefer~\cite{Sch78}.
\begin{lemma}\label{lemma: in pspace}
     Both Client-Waiter and Waiter-Client positional games are in \PSPACE{}.
\end{lemma}

\begin{proof} 
Let $H = (V,E)$ be a hypergraph. Each turn, Waiter has the choice among at most $\binom{|V|}{2}$ moves and Client has the choice among $2$ moves. Since the game ends in at most $|V|$ moves, it is in \PSPACE\  using the same proof as in Lemma~2.2 from Schaefer~\cite{Sch78}.
\end{proof}

\subsection{Quantified Boolean Formula and paired SAT}
\label{sec:PairedSAT}

The most classical \PSPACE-complete problem is 3-\QBF, the quantified version of \SAT. Our hardness proof is a reduction from \(3\)-\QBF, but since the roles of the players in Client-Waiter games are very different, we introduce an intermediate problem \PS.

First we recall the definition of 3-\QBF, in its gaming version, as it was done by Rahman and Watson~\cite{RW21}, and later Gledel and Oijid~\cite{GO23} to prove that Maker-Breaker and Avoider-Enforcer games are \PSPACE-hard respectively.

Given a 3-CNF quantified formula $\varphi = \exists x_1, \forall y_1, \dots, \exists x_{n} \forall y_{n} \psi$, where $\psi$ is a $3$-CNF without quantifier, the 3-\QBF game is played by two players, Satisfier and Falsifier. Satisfier chooses the value of $x_1$, then Falsifier chooses the value of $y_1$, and so on until the last variable has its value chosen. At the end, a valuation $\nu$ of the variables is obtained, and Satisfier wins if and only if $\nu$ satisfies $\psi$.

\begin{theorem}[Stockmeyer and Meyer~\cite{SM73}]\label{QBF PSPACE COMPLETE}
Determining if Satisfier has a winning strategy in the $3$-\QBF game is \PSPACE-complete.
\end{theorem}

The game \PS is introduced in the introduction (Definition~\ref{def: paired-sat}). This is a variant of \(3\)-\QBF where, at each turn, Satisfier chooses an index \(i\) and instantiates the variable \(x_i\), and then Falsifier instantiates the variable \(y_i\). The main idea of this game is to introduce a CNF-game mimicking the fact that one player chooses which variable the second player has to play on.


\begin{theorem}
    Determining the winner of the \PS-game is \PSPACE-complete.
\end{theorem}

\begin{proof}
Since at any moment of the game, the number of options for a player is at most $n$ and the game runs for at most $n$ turns, determining the winner of the \PS-game is in \PSPACE\ using the same proof as Lemma~2.2 from Schaefer~\cite{Sch78}.

We provide a reduction from 3-\QBF. Let $\psi = \exists x_1 \forall y_1 \cdots \exists x_n \forall y_n, \varphi$ be a \QBF formula. We construct an instance of the \PS-game $(\varphi', X)$ as follows: 

\begin{itemize}
    \item $X = \{(z_0, y_0), (x_1, t_1), (z_1, y_1), \dots, (x_n, t_n), (z_n, y_n)\}$, where $y_0$, the $z_i$s, and the $t_j$s are new variables.
    \item $\varphi' = \varphi \wedge \bigwedge\limits_{1 \le i \le n} (y_{i-1} \oplus t_i \oplus z_i)$. 
    
    where $\oplus$ refers to the XOR operator:
    \[a \oplus b \oplus c = (a \vee b \vee c) \wedge (\neg a \vee \neg b \vee c) \wedge (\neg a \vee b \vee \neg c) \wedge (a \vee \neg b \vee \neg c).\]
\end{itemize}

We prove that Satisfier wins on $\psi$ if and only if she wins on $(\varphi', X)$.

First note that whenever two values of $a \oplus b \oplus c$ are known, the player who chooses the last value can always decide to satisfy or not $(a \oplus b \oplus c)$.
Suppose first that Satisfier has a winning strategy $\strat$ on $\psi$, and consider the following strategy for him on $(\varphi', X)$:

\begin{itemize}
    \item Satisfier instantiates the pairs $(z_0, y_0), (x_1, t_1), (x_2, t_2), \dots, (z_n, y_n)$ in that order.
    \item Whenever Satisfier has to choose a value for a variable $x_i$, he follows $\strat$ with the corresponding values of the $x_j$ and $y_j$ for $1 \le j < i$. This is always possible as the order was given above.
    \item Whenever Satisfier has to give a value to a variable $z_i$, she gives the value so that $y_{i-1} \oplus t_i \oplus z_i$ is satisfied (if $i=0$ she can instantiate the variable $z_0$ by either $\top$ or $\bot$).
\end{itemize}

Following this strategy, all the clauses $(y_{i-1} \oplus t_i \oplus z_i)$ are satisfied as Satisfier always chooses the last vertex of these clauses (which appears in exactly one of them), and as $\strat$ has a winning strategy in $\psi$, it satisfies $\varphi$, as the variables are chosen in the same order.

Now suppose that Falsifier has a winning strategy $\strat$ on $\psi$, and consider the following strategy on $(\varphi', X, Y)$:

\begin{itemize}
    \item While Satisfier plays the pairs following the order $(z_0, y_0), (x_1, t_1), (x_2, t_2), \dots, (z_n, y_n)$, Falsifier gives to the corresponding $t_i$ the value $\bot$, and to $y_i$ the value given by $\strat$.
    \item If Satisfier plays a pair $(x_j, t_j)$ before she should, Falsifier still gives the valuation $\bot$ to $t_j$ and ignores this move while choosing the valuations of the $y_i$ for $i\le j$.
    \item If Satisfier plays a vertex $z_j$ before he should, the first time it happens, all the unplayed variables in the clause $(y_{j-1} \oplus t_j \oplus z_j)$ will be played by Falsifier. Therefore, Falsifier can just win by choosing a good valuation for $y_{j-1}$ and $t_j$.
\end{itemize}

Following this strategy, if Satisfier instantiates a variable $z_i$ whereas there is $j\le i$ such that the variable $x_j$ has not yet been instantiated then Falsifier wins. Indeed, the first time it happens, either $z_{i-1}$ or $x_i$ has not been instantiated (otherwise it already happened when Satisfier instantiated $z_{i-1}$). Consequently, Falsifier wins through the clause $(y_{i-1} \oplus t_i \oplus z_i)$. Otherwise, each time Falsifier has to choose a valuation for a variable $y_i$, all the vertices $x_j$ with $j \le i$ have already been played and so, he can play according to $\strat$. As $\strat$ is a winning strategy, in both case, Falsifier can make a clause unsatisfied and therefore wins the game.
\end{proof}

\subsection{Reduction to Client-Waiter games}

\subsubsection{Blocks in Client-Waiter games}\label{subsection blocks}

\begin{figure}[ht]
    \centering

\begin{tikzpicture}

\draw (-1,.5) node[v] (v5)  {};
\draw (-2,.5) node[v] (v6)  {};

\draw (0,0) node[v] (v1)  {};
\draw (0,1) node[v] (v2)  {};
\draw (1,0) node[v] (v3)  {};
\draw (1,1) node[v] (v4)  {};

\draw[color = red] \hedgeiii{v1}{v2}{v3}{4.5mm};
\draw[color = blue] \hedgeiii{v2}{v4}{v3}{5.25mm};
\draw[color = green] \hedgeiii{v4}{v3}{v1}{3.75mm};
\draw[] \hedgeiii{v1}{v2}{v4}{3mm};
\draw[] \hedgeii{v5}{v6}{4mm};

\draw[] \hedgeiii{v5}{v2}{v1}{4mm};

\end{tikzpicture}

    \caption{A blocks-hypergraph. The two vertices on the left form a block. The four on the right a second one. The hyperedge between them is in no block}
    \label{fig:block}
\end{figure}
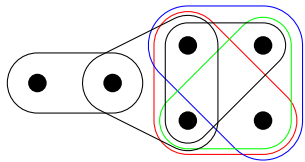

The main tool of several reductions of positional games is pairing strategies. However, this cannot be applied to Client-Waiter games, since only Waiter has choices about how to make the pairs. We present {\em blocks-hypergraphs} and {\em block-strategies} that will be used similarly to pairing strategies to ensure that client can claim some vertices. A blocks-hypergraph is depicted in Figure~\ref{fig:block}. 
The idea of blocks was already used in the \NP-hardness proof from \cite{CMP11}, but we present here a more formal definition of them. Intuitively, Blocks are a generalisation of Lemma~\ref{lem:size2hyperedge}, which are blocks of size $2$.

\begin{definition}[Blocks] \label{definition: blocks in Client-Waiter}
Let $\hxf$ be a hypergraph. A {\em block} $B \subset \som$ of size $2k$ is a set of vertices such that $|B| = 2k$ for some $k \ge 1$, and any set of $k+1$ vertices of $B$ is an \hedge.

If $\hyp$ can be partitioned into blocks, we say that $\hyp$ is a {\em blocks-hypergraph}.
\end{definition}\index{blocks-hypergraph} \index{Block}

\begin{lemma}
    Let $\hxf$ be a hypergraph, and let $B$ be a block of $\hyp$. If Waiter has a winning strategy in $\hyp$, she has to offer the vertices of $B$ two by two.
\end{lemma}

\begin{proof}
    Suppose that Waiter has a winning strategy in which she does not offer all vertices of $B$ two by two. Let $k = \frac{|B|}{2}$. The first time she presents a vertex $x\in B$ with a vertex $y \notin B$, Client can choose $x$. Then, each time Waiter offers at least one vertex in $B$, Client claims it. In the end, Client will claim at least $k+1$ vertices of $B$ and therefore wins. Thus, if Waiter has a winning strategy, she has to offer the vertices of $B$ two by two.
\end{proof}

\begin{corollary}\label{corollary: always play in blocks}
    Let $\hxf$ be a blocks-hypergraph. If Waiter has a winning strategy in $\hyp$, any pair of vertices she offers belong to a same block of $\hyp$.
\end{corollary}

\subsubsection{Construction of the hypergraph} \label{subsection construction client-waiter}

We show now the reduction. The main idea of the reduction is that we construct a blocks-hypergraph, such that each block corresponds to the valuation that will be given to a variable.

Let $(\varphi,X)$ be an instance of \PS where $X = \{(x_1, y_1),\ldots,(x_n, y_n)\}$, and $\varphi = \bigwedge_{1 \le j \le m} C_j$ is a 3-CNF on the variables of $X$. We build a hypergraph $\hxf$ as follows.

Let us define the set $\som$ of $8n$ vertices. Let $\som = \bigcup_{1 \le i \le n} S_i \cup F_i $, with for $1 \le i \le n$, $S_i = \{s_{i}^0,s_{i}^T,s_{i}^F,s_{i}^1\}$ (gadget which encodes Satisfier's choice for the variable $x_i$) and $F_i = \{f_{i}^0,f_{i}^T,f_{i}^{T'},f_{i}^F\}$ (gadget which encodes Falsifier's choice for the variable $y_i$).

Now we focus on the construction of the \hedges. 

\begin{figure}
    \centering

\begin{tikzpicture}

    \draw (0,0) node[v] (si0)  {} node[above= .15cm] {$s_i^0$};
    \draw (-1,-1) node[v] (sit)  {} node[above= .15cm] {$s_i^T$};
    \draw (1,-1) node[v] (sif)  {} node[above= .15cm] {$s_i^F$};
    \draw (0,-2) node[v] (si1)  {} node[above= .15cm] {$s_i^1$};

    \draw (-1,2) node[v] (fit)  {} node[above= .15cm] {$f_i^T$};
    \draw (1,2) node[v] (fitp)  {} node[above= .15cm] {$f_i^{T'}$};
    \draw (0,3) node[v] (fi0)  {} node[above= .15cm] {$f_i^0$};
    \draw (2,1) node[v] (fif)  {} node[above= .15cm] {$f_i^F$};

\draw[style = dashed] \hedgeiiii{si0}{sif}{si1}{sit}{5mm} node[below right = .2cm] {Block $S_i$};
\draw[style = dashed] \hedgeiiii{fi0}{fitp}{fif}{fit}{5mm} node[right = .2cm]{Block $F_i$};
\draw[color = red] \hedgeiiii{fit}{fi0}{si0}{sit}{4mm};
\draw[color = red] \hedgeiiii{fitp}{sif}{si0}{fi0}{3mm};
\draw[color = blue] \hedgeiii{fif}{sit}{fitp}{3mm};
\draw[color = blue] \hedgeiiii{fif}{sif}{si0}{fit}{4mm};

\end{tikzpicture}

    \caption{Gadget for the vertices in $B_i$. A dashed set represents a block, i.e.\ all hyperedges of size three are present in it.}
    \label{fig: gadget client waiter}
\end{figure}
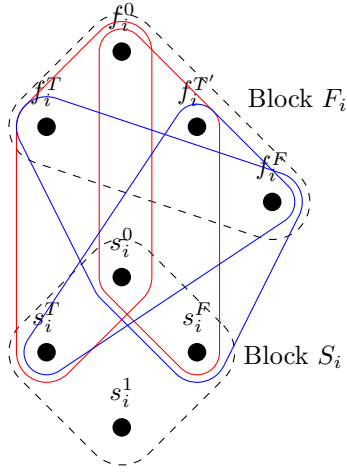

\begin{itemize}[noitemsep]
    \item The block-\hedges $B = \bigcup_{1 \le i \le n} B_i$, which make each $S_i$ and each $F_i$ a block:
\[
    B_i = \left\{ H \subseteq S_i \mid  \lvert H \rvert=3\} \cup \{H \subseteq F_i \mid  \lvert H \rvert=3 \right\}.
\]

\item The pair-\hedges $P = \bigcup_{1 \le i \le n} P_i$ (see Figure~\ref{fig: gadget client waiter}):
\[ 
    P_i = \left\{ \{s_i^0, s_i^T, f_i^0, f_i^T \}, \{s_i^0, f_i^F, f_i^T, s_i^{F} \}, \{s_i^0, f_i^F, s_i^T, f_i^{T'} \}, \{s_i^0, s_i^F, f_i^0, f_i^{T'} \} \right \}.
\]

\item The clause-\hedges. Each clause $C_j \in \varphi$ is a set of three literals $\{\ell_{j}^1,\ell_{j}^2,\ell_{j}^3\}$. We define first, for $1 \le j \le m$ and $k \in \{ 1, 2, 3\}$, the set $H_{j}^k$ which encodes the property that the literal $\ell_{j}^k$ is instantiated to $\bot$. 
\begin{align*}
    H_{j}^k = \begin{cases} 
    \{\{s_{i}^0,s_{i}^T\}\} & \text{if } \ell_{j}^k = x_i \\
    \{\{s_{i}^0,s_{i}^F\}\} & \text{if } \ell_{j}^k = \neg x_i \\
    \{\{f_{i}^0,f_{i}^T\},\{f_{i}^0,f_{i}^{T'}\}\} & \text{if } \ell_{j}^k =  y_i \\
    \{\{f_{i}^F\}\} & \text{if } \ell_{j}^k = \neg y_i. 
    \end{cases}
\end{align*}
We define now the set of \hedges:
   
$$ C = \bigcup_{C_j \in \varphi} H_j. $$

\noindent with $H_j = \left\{ h_{1} \cup h_{2} \cup h_{3} \mid \forall k \in \{ 1,2,3 \},\, h_k \in H_{j}^k\right\}$

For example, if $C_j = x_1 \vee y_1 \vee \neg y_2$, we have $H_j^1 = \{ \{s_1^0, s_1^T\} \}$, $H_j^2 = \{ \{f_1^0, f_1^T\}, \{ f_1^0, f_1^{T'}\} \}$ and $H_j^3 = \{\{ f_2^F\}\}$. Finally, we have two \hedges to encode $C_j$: $H_j = \{ (s_1^0, s_1^T, f_1^0, f_1^{T}, f_2^F), (s_1^0, s_1^T, f_1^0, f_1^{T'}, f_2^F)  \}$

\end{itemize}

The gadget for the pair $(x_i,y_i)$ is depicted in Figure~\ref{fig: gadget client waiter}.

Intuitively, these \hedges are constructed in such a way that:
\begin{itemize}
    \item The block-\hedges $B$ force Waiter to always propose two vertices in the same set. 
    \item The pair-\hedges $P$ force Waiter to give to Client the choice of the value of $y_i$ after she has made her choice for $x_i$.
    \item The clause-\hedges $C$ which represent the clauses of $\varphi$ and make the equivalence between a win of Waiter and a valuation that satisfies $\varphi$.
\end{itemize}

Finally, the reduction associates to the instance $(\varphi,X)$ of {\sc Paired SAT} the hypergraph $\hyp = (\som, \WS$) with $\WS = B \cup P \cup C$ of Client-Waiter. This reduction is polynomial as $B$ contains $8n$ \hedges, $P$ contains $4n$ \hedges, and $C$ contains at most $8m$ \hedges (where $m$ is the number of clause of $\varphi$).

We define an underlying assignment of the variables, corresponding to the moves on the hypergraph as follows:
\begin{itemize}[noitemsep]
    \item If Client claims $s_i^T$, $x_i = \bot$
    \item If Client claims $s_i^F$, $x_i = \top$
    \item If Client claims $f_i^0$ and one of $f_i^T$, $f_i^{T'}$, $y_i = \bot$
    \item If Client claims $f_i^F$, $y_i = \top$
\end{itemize}

We prove in next sections (Lemmas~\ref{falsifier implies client} and~\ref{falsifier implies client}) that Waiter has a winning strategy on \(\hyp\) if and only if Satisfier has a winning strategy on \(\varphi\). Then, we can obtain the proof of Proposition~\ref{prop: Client-Waiter pspace complete}.

\begin{proof}[Proof of Proposition~\ref{prop: Client-Waiter pspace complete}]

First, according to Lemma~\ref{lemma: in pspace}, Client-Waiter games are in \PSPACE. We prove the hardness by a reduction from \PS.

Let $(\varphi,X)$ be an instance of \PS. Consider the hypergraph $\hyp$ obtained from the reduction provided in Subsection~\ref{subsection construction client-waiter}. It has $O(|X|)$ vertices and $O(|X| + |\varphi|)$ \hedges, which is polynomial. According to Lemma~\ref{Satisfier implies waiter} and Lemma~\ref{falsifier implies client}, Satisfier wins in $(\varphi, X)$ if and only if Waiter wins in $\hyp$. Therefore, determining the winner of a Client-Waiter game is \PSPACE-complete.

Moreover, as any \hedge of $\hyp$ has size at most~$6$, the problem is even \PSPACE-complete restricted to hypergraphs of rank~$6$.
\end{proof}

\subsubsection{Waiter's winning strategy} \label{subsection waiter}

In this section, we prove that if Satisfier has a winning strategy in $\varphi$ for the \PS-game, then Waiter has a winning strategy in $\hyp$ for the Client-Waiter game. 

\begin{lemma}\label{Satisfier implies waiter}
    If Satisfier has a winning strategy in $\varphi$, then Waiter has a winning strategy in $\hyp$. 
\end{lemma}

\begin{proof}
    Let $\strat$ be a winning strategy for Satisfier, consider a strategy for Waiter as follows. If $\strat$ selects an integer $1 \le i \le n$, and puts $x_i$ to $\top$ (resp. $\bot$), Waiter plays in the block \(B_i\) and selects the pair $(s_i^0, s_i^T)$ (resp. $(s_i^0, s_i^F)$). Then, she plays the pair corresponding to the two other vertices in the block $S_i$. To determine the value of $y_i$, she plays $(f_i^F,f_i^{T})$ (resp. $(f_i^F, f_i^{T'}) $) and finally, the remaining pair of the block $F_i$. If Client chooses $f_i^F$, she considers that $y_i = \top$, otherwise, she considers that $y_i = \bot$ in $\strat$.

As this strategy always propose vertices in blocks, Client cannot win with the \hedges in $B$:
\begin{itemize}
    \item In the case Waiter offers $(s_i^0, s_i^T)$, if Client does not choose $s_i^0$, as it is in all the \hedges of $P_i$, Waiter can not lose on $P_i$. Otherwise, Waiter claims $s_i^T$, and proposes the pairs $(f_i^F, f_i^T)$ and $(f_i^{T'}, f_i^0)$ which cover all the remaining \hedges of $P_i$. 
    \item In the other case Waiter offers $(s_i^0, s_i^F)$, again if Client chooses $s_i^0$, Waiter cannot lose on $P_i$. Otherwise, Waiter claims $s_i^F$, and proposes the pairs $(f_i^F, f_i^{T'})$ and $(f_i^T, f_i^0)$  which covers all the remaining \hedges of $P_i$. 
\end{itemize}

We now consider the clause-\hedges:

Let $C_j$ be a clause of $\varphi$. It is sufficient to prove that there exists $1 \le k \le 3$ such that Waiter claims a vertex in $H_j^k$. As $\strat$ is a winning strategy for Satisfier in $\varphi$, there exists, at the end of the game, an index \(1 \le k \le 3\) such that the assignment of the literal $\ell_j^k$ satisfies $C_j$. 

\begin{itemize}
    \item If $\ell_j^k = x_i$, by construction of the strategy, Waiter offers the pair $(s_i^0, s_i^T)$, therefore she claims a vertex in $H_j^k$.
    \item If $\ell_j^k = \neg x_i$, by construction of the strategy, Waiter offers the pair $(s_i^0, s_i^F)$, therefore she claims a vertex in $H_j^k$.
    \item If $\ell_j^k = y_i$, by construction of the strategy, Waiter has considered that $y_i$ was put to $\top$ according to the choices of Client, which corresponds to the case where Client has chosen $f_i^F$. Therefore, she has claimed two of the three vertices $\{f_i^0, f_i^T, f_i^{T'}\}$, and so has a vertex in each set of $H_j^k$. 
    \item If $\ell_j^k = \neg y_i$, by construction of the strategy, Waiter has considered that $y_i$ was put to $\bot$ according to the choice of Client, which corresponds to the case where Client has not chosen $f_i^F$, therefore Waiter has claimed it, and thus has a vertex in $H_j^k$.
\end{itemize}

Finally, Client can not fill up an \hedge in $H_j$ for any $1 \le j \le m$. So the described startegy is winning for Waiter.
\end{proof}

\subsubsection{Client's winning strategy} \label{subsection client}

We prove now the other direction.

\begin{lemma}\label{falsifier implies client}
If Falsifier has a winning strategy in $\varphi$, then Client has a winning strategy in $\hyp$.
\end{lemma}

\begin{proof}
Suppose now that Falsifier has a winning strategy $\strat$ in $\varphi$. We provide a strategy for Client. 

First, notice that, as $\hyp$ is a blocks-hypergraph, we can suppose that Waiter always offer vertices in blocks, according to Corollary~\ref{corollary: always play in blocks}. Moreover, since each block has size $4$, once the first pair of a block is offered, Client knows that the second pair will be proposed at some point and can already decide her move on this second pair. Therefore it is sufficient to have a winning strategy when Waiter always offers simultaneously the two pairs of a same block.

Let $1 \le i \le n$ be an integer, and suppose that Waiter offers two disjoint pairs of a block $S_i$ or $F_i$. 

\begin{itemize}
    \item If the pairs are in $S_i$, Client can ensure to claim $s_i^0$ and one of $s_i^T$ and $s_i^F$. If he has claimed $s_i^F$, he considers that Satisfier has instantiated $x_i$ to $\top$ and if he has claimed $s_i^T$, he considers \(x_i\) has been instantiated to $\bot$. 
    \begin{itemize}
        \item Assume Client claims \(s_i^0\) and \(s_i^T\). Since \(\{s_i^0,s_i^T,f_i^T,f_i^0\}\) is a \hedge of \(P_i\), by Lemma~\ref{lem:size2hyperedge} Waiter will offer the pairs \(\{f_i^T,f_i^0\}\) and \(\{f_i^{T'},f_i^F\}\). 
        Client can follow $\strat$ by picking \(f_i^0\) and \(f_i^{T'}\) if $\strat$ instantiates \(y_i\) to \(\bot\), and by picking \(f_i^T\) and \(f_i^F\) otherwise.  
        \item If Client claims \(s_i^0\) and \(s_i^{F}\), the result is similar by switching the roles of \(f_i^T\) and \(f_i^{T'}\). 
    \end{itemize}  
    \item If the pairs are in $F_i$, as $F_i$ is a block of four vertices, Waiter has only three possible ways to pair them.
    \begin{itemize}
        \item If the pairs are $\{f_i^0, f_i^T\}$ and $\{f_i^F, f_i^{T'}\}$, Client will claim either the two vertices \((f_i^{T},f_i^{F})\), or the two vertices \((f_i^0,f_i^{T'})\). In both cases, Waiter will have to offer the pairs \(\{s_i^0,s_i^F\}\) and \(\{s_i^{T},s_i^1\}\), and Client will claim \(s_i^0\) and \(s_i^F\). This case is already described above, and so the choice of the claiming pair \((f_i^{T},f_i^{F})\) or \((f_i^0,f_i^{T'})\) can be done as previously in following \(\strat\).
        \item If the pairs are $\{f_i^T, f_i^{T'}\}$ and $\{f_i^0, f_i^F\}$, Client can still claim either \((f_i^{T},f_i^{F})\), or \((f_i^0,f_i^{T'})\) and the case is identical.
        \item Otherwise the pairs are $\{f_i^T, f_i^F\}$ and $\{f_i^0, f_i^{T'}\}$. Client will claim either \((f_i^{T'},f_i^{F})\), or \((f_i^0,f_i^{T})\). In both cases, Waiter will have to offer the pairs \(\{s_i^0,s_i^T\}\) and \(\{s_i^{F},s_i^1\}\), and Client will claim \(s_i^0\) and \(s_i^T\). This is again a case described above.
    \end{itemize}
\end{itemize}

Following this strategy until the end, the underlying valuation of the claimed vertices is the one obtained by $\strat$ in $\varphi$. By hypothesis, there exists a clause $C_j \in \varphi$ in which all literals are set to $\bot$. 
Thus Client wins by filling up the \hedge \(H_j\).
\end{proof}

\subsection{Reduction to 6-uniform hypergraphs} \label{subsection conclusion}

As it is done by Rahman and Watson for Maker-Breaker games~\cite{RW21}, or by Gledel and Oijid for Avoider-Enforcer games~\cite{GO23}, Proposition~\ref{prop: Client-Waiter pspace complete} can be strengthen by transforming the hypergraph of rank \(6\) into a \(k\)-uniform one (with \(k \ge 6\)). We achieved this with the following lemma:

\begin{lemma}\label{to 6uniform Client-Waiter}
    Let $\hxf$ be a hypergraph of rank $k$. Let $m = \min_{e\in E} (|e|)$. If $m <k$, there exists a hypergraph $\hyp'= (\som',\WS')$ of rank~$k$ where $\min_{e\in E} (|e|) = m+1$, having $|\WS'|\le |\WS| + \binom{2k-2}{k}$ and $|\som'| \le |\som| + 2(k-1)$ such that Client has a winning strategy in the Client-Waiter game on $\hyp$ if and only if he has one in $\hyp'$.
\end{lemma}

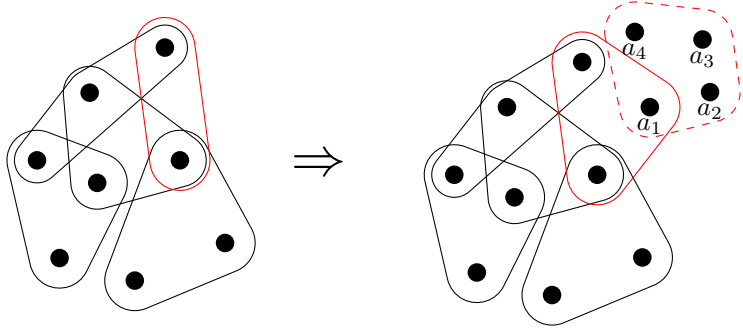
\begin{figure}[ht]
    \centering

\begin{tikzpicture}
    
    \node at (0,0){
        \begin{tikzpicture}
        
            \node[v] (a) at (0,0){};
            \node[v] (b) at (0.5,1){};
            \node[v] (c) at (-0.3,1.3){};
            \node[v] (d) at (1,-0.3){};
            \node[v] (e) at (0.4,2.2){};
            \node[v] (f) at (1.4,2.8){};
            \node[v] (g) at (1.6,1.3){};
            \node[v] (h) at (2.2,0.2){};
            
            \draw \hedgeiii{a}{c}{b}{4mm};
            \draw \hedgeiii{c}{e}{f}{3mm};
            \draw \hedgeiii{b}{e}{g}{3.5mm};
            \draw \hedgeiii{d}{g}{h}{4mm};
            \draw[color=red] \hedgeii{g}{f}{4mm};
            
        \end{tikzpicture}

    };

    \node at (2.5,0){\huge $\Rightarrow$};

    \node at (6,0){
    
        \begin{tikzpicture}
        
            \node[v] (a) at (0,0){};
            \node[v] (b) at (0.5,1){};
            \node[v] (c) at (-0.3,1.3){};
            \node[v] (d) at (1,-0.3){};
            \node[v] (e) at (0.4,2.2){};
            \node[v] (f) at (1.4,2.8){};
            \node[v] (g) at (1.6,1.3){};
            \node[v] (h) at (2.2,0.2){};
            
            \node[v] (u1) at (2.1,3.2){};
            \node[v] (u2) at (2.3,2.2){};
            \node[v] (u3) at (3.1,2.4){};
            \node[v] (u4) at (3,3.1){};
            
            \draw \hedgeiii{a}{c}{b}{4mm};
            \draw \hedgeiii{c}{e}{f}{3mm};
            \draw \hedgeiii{b}{e}{g}{3.5mm};
            \draw \hedgeiii{d}{g}{h}{4mm};
            \draw[color=red] \hedgeiii{g}{f}{u2}{4mm};
            \draw[color=red, dashed] \hedgeiiii{u2}{u1}{u4}{u3}{4mm};
            
            \node[below=1pt] at (u2){$a_1$};
            \node[below=1pt] at (u3){$a_2$};
            \node[below=1pt] at (u4){$a_3$};
            \node[below=1pt] at (u1){$a_4$};

        \end{tikzpicture}

    };


\end{tikzpicture}

    \caption{The construction of Lemma 22 with $k= 3$ and $m= 2$.  The dashed set is a block and contains the four hyperedges of size 3.  The resulting hypergraph is 3-uniform}
    \label{fig:to6unif CW}
\end{figure}

\begin{proof}
    Let $\hxf$ be a hypergraph of rank~$k$. We construct the hypergraph $\hyp' = (\som', \WS')$. Let $A = \{a_1, \dots, a_{2(k-1)}\}$ be $2(k-1)$ new vertices, and set $\som' = \som \cup A$. We make $A$ a block, i.e.\ we define the set $ U = \{B \subset A \mid \lvert B\rvert = k \}$ of \hedges. Then, we introduce $L = \{e \cup \{a_1\} \mid e\in \WS \text{ and } \lvert e\rvert = m  \}$. Finally, we define our \hedges as follows:

\[ 
\WS' = \{e\in \WS \mid \lvert e \rvert \ge m+1 \} \cup L \cup U.
\]

The construction is depicted in Figure~\ref{fig:to6unif CW}.

Client wins in $\hyp$ if and only if he wins in $\hyp'$. Indeed, suppose that Client wins in $\hyp$, as $A$ is a block, Waiter will have to offer the vertices of $A$ two by two. Therefore, Client can claim $a_1$ when it will be proposed and apply the strategy for $\hyp$ on vertices of \(V\). Then, if Client fills up an \hedge $e \in \WS$ with $|e| = m$, he also fills up $e \cup \{a_1\}$ in $\hyp'$, and if she fills up \(e \in \WS\) with $|e| \ge m+1$, she also fills up $e$ in $\hyp'$. Reciprocally, if Waiter has a winning strategy in $\hyp$, she can offer the same pairs in $\hyp'$ and the vertices of $A$ two by two. Then, for any \hedge $e \in \WS'$, if $e \in \WS$ or $e \setminus \{a_1\} \in \WS$, Waiter has claimed a vertex of it by her winning strategy in $\hyp$. Otherwise, $e \in U$, and she claims one vertex of it since she claims $k-1$ vertices of $A$.
\end{proof}

Hence, we obtain Theorem~\ref{6-uniform Client-Waiter} by combining Proposition~\ref{prop: Client-Waiter pspace complete} and Lemma~\ref{to 6uniform Client-Waiter}.

\begin{proof}
    The hypergraph obtained in the proof of Proposition~\ref{prop: Client-Waiter pspace complete} has rank~6. 
    Therefore, by applying \(k-1\) times Lemma~\ref{to 6uniform Client-Waiter}, with $m = 1, \ldots, k-1$, we obtain a \(k\)-uniform hypergraph having at most $2(k-1)^2$ more vertices and $(k-1)\binom{2k-2}{k}$ more \hedges. Thus, when \(k\) is fixed, this construction is still polynomial, and the hypergraph obtained is $k$-uniform. 
\end{proof}

\section{Rank 3 Client-Waiter equivalent to the problem of detecting a tadpole} \label{sec: rank3 client-waiter}

Similarly to Maker-Breaker, the \PSPACE-hardness of Client-Waiter games restricted to $6$-uniform hypergraphs leads us to consider smaller rank hypergraphs. In this section, we first provide a linear time algorithm to compute the winner in rank~2 hypergraphs, then we prove that rank~3 hypergraphs are in $\Delta_2^P$.

\subsection{Rank~2 games}

Let us start with the easier case of rank~2 hypergraphs. We prove that Client-Waiter games are tractable restricted to them.

\begin{proposition}\label{prop: rank 2 Client-Waiter}
        Let $\hxf$ be a rank~$2$ hypergraph. The winner of a Client-Waiter game played on $\hyp$ can be computed in linear time.
\end{proposition}

\begin{proof}
    We prove that Client wins on $\hyp$ if and only if there exists an \hedge of size $1$, or if two \hedges intersect. The result directly follows since these properties can be checked in linear time.

    \begin{itemize}
        \item  If $\hyp$ has an \hedge $\{x\}$ of size $1$, when Waiter proposes the vertex $x$, Client can claim it and win. As the last vertex goes to Client if the number of vertices if odd, Client can ensure to get \(x\) and win.
        \item If there exist two \hedges $\{a,b\}$ and $\{a,c\}$ which intersect, Client applies the following strategy: \begin{itemize}[noitemsep]
            \item Whenever Waiter offers two of these vertices, Client claims one of them, and he takes $a$ if it is available.
            \item If Waiter offers one of these vertices with any other, Client claims the one in $\{a,b,c\}$. 
            \item Otherwise, Client takes any vertex.
        \end{itemize}
        Following this strategy, Client claims \(a\) and at least one vertex from $\{b,c\}$. Therefore, he wins. Once again, even if the number of vertices is odd, this strategy can be applied.
    \end{itemize}

    Reciprocally, if all \hedges have size $2$ and do not intersect (\(H\) is in fact a graph, and even a matching), by always selecting a pair $(a,b)$ such that $(a,b)$ is an \hedge, Waiter gets one vertex from each \hedge and wins.
\end{proof}

\subsection{Rank~3 games}

We now focus on Client-Waiter games restricted to rank~$3$ hypergraphs. We show that testing if a position is winning for Client in a hypergraph \(\hyp\) of rank \(3\) reduces to the problem of searching two structures (the \(a\)-snakes and the \(ab\)-tadpoles) in \(\hyp\). We start by defining them.




\begin{definition}\label{def_path_cycles}
    Let \(\hyp=(V,E)\) be a hypergraph and \(a, b \in V\). 
    
    A sequence of \hedges of \(\hyp\), \(\mathcal{P} = (e_1,\ldots,e_t)\), is an \(ab\)-{\em path} if \(a \in e_1\), \(b \in e_t\), and for \(1 \le i \le t-1\) the \hedges \(e_i\) and \(e_{i+1}\) intersect. The number of \hedges \(t\) is called the {\em length} of \(\mathcal{P}\). An \(a\)-{\em cycle} is an \(aa\)-path of length at least two. An \(ab\)-path is said {\em linear} if the size of the intersection of two consecutive \hedges is always one. An \(a\)-cycle of length at least \(3\) is {\em linear} if the \(aa\)-path is linear and if the intersection of the end \hedges is exactly \(\{a\}\). We continue to call {\em linear} an \(a\)-cycle \((e_1,e_2)\) of length \(2\) if \(\lvert e_1 \cap e_2\rvert = 2\). An \(ab\)-path is said {\em simple} if \(a\) appears only in \(e_1\), \(b\) appears only in \(e_t\), and if whenever \(e_i\) and \(e_j\) intersect, then \(\lvert i-j\rvert \le 1\). Similarly, an \(a\)-cycle is simple if whenever \(e_i\) and \(e_j\) intersect, then \(\lvert i-j\rvert \le 1\) or \(\{i,j\}=\{1,t\}\).
    
    An \(ab\)-path is also called an \(a\)-path or a path. A cycle is \(a\)-cycle for some \(a\) in \(V\). 
\end{definition}

\begin{definition}\label{def_snake}
    Let \(a\) be a vertex of \(\hyp\). An \(a\)-{\em snake} is an \(a\)-path \((e_1,\ldots,e_t)\) with \(t \ge 1\) such that for all \(1 \le i \le t-1\), $e_i$ has exactly size \(3\), and \(e_t\) has size at most \(2\).
\end{definition}

\begin{definition}\label{def_tadpole}
    Let \(H=(V,E)\) be a hypergraph and \(a, b \in V\). 
	If \(a\neq b\), an $ab$-{\em tadpole} is a sequence of \hedges $T=(e_1,\ldots,e_s, f_1, \ldots, f_t)$ where:
	\begin{itemize}[noitemsep,nolistsep]
		\item $a$ belongs to \(e_1\) and no other \hedge;
		\item $b$ belongs to \(e_s\), \(f_1\), \(f_t\) and no other \hedge;
		\item $(e_1,\ldots,e_s)$ is a \(3\)-uniform simple linear $ab$-path $\mathcal{P}_T$;
		\item $(f_1, \ldots, f_t)$ is a \(3\)-uniform simple linear $b$-cycle $\mathcal{C}_T$;
		\item \(b\) is the only vertex which appears both in $\mathcal{P}_T$ and $\mathcal{C}_T$.
	\end{itemize}
	If \(a=b\), an \(ab\)-tadpole is just a \(3\)-uniform simple linear \(a\)-cycle.
	When \(T\) is an \(ab\)-tadpole, we may simply say $T$ is an \(a\)-\textit{tadpole}, or even just a \textit{tadpole}.
\end{definition}

We consider the problem \textsc{Tadpole}: Given a vertex \(u\) in a \(3\)-uniform hypergraph \(H\), decide if there is a \(u\)-tadpole in \(H\).

We denote by \(o(H)\) the outcome of an optimal Client-Waiter game on the hypergraph \(H\). We will write \(o(H)= \Waiter\) when Waiter has a winning strategy on \(H\), and \(o(H) = \Client\) otherwise.

Let \(u\) be a vertex of \(H\), we consider the following family \(\uCwin_H\) of subhypergraphs of \(H\):
\[
	T \in \uCwin_H \iff \begin{cases}
			T \text{ is a } u\text{-snake} \\
			\text{or } T \text{ is a } u\text{-tadpole}.
		\end{cases}
\]

Notice that \(3\)-uniform simple linear \(u\)-cycles are particular cases of \(u\)-tadpoles, so such subhypergraphs are also in \(\uCwin_H\). When the hypergraph \(H\) is known, we will simply write \(\uCwin\).

We notice that Waiter has a winning strategy in \(H=(V,E)\) which starts by offering a pair \(\{u,v\}\) if and only if Waiter has a winning strategy in both trace hypergraphs \(H_1 = T_{V \setminus \{u,v\}}(H_{| V \setminus \{u\}})\) and \(H_2 = T_{V \setminus \{u,v\}}(H_{| V \setminus \{v\}})\). So to simplify notations, we will write \(H^{+v}\) for the trace \(T_{V\setminus\{v\}}(H)\) and \(H^{-v}\) for the induced subhypergraph \(H_{|V\setminus \{v\}}\). In particular, we can just write \(H_1 = H^{+v-u}\) and \(H_2 = H^{+u-v}\).

\begin{lemma}\label{lem:gTadpoleCwin}
	Let \(u\) be a vertex of \(H\) hypergraph of rank \(3\) and let \(T \in \uCwin\). Then \(o(T^{+u}) = \Client\).
\end{lemma}

\begin{proof}
	We do the proof by induction on the size (number of \hedges) of \(T\). By definition, each subhypergraph of \(\uCwin\) contains at least one \hedge. 
	
	Assume that \(T \) contains exactly one \hedge. It means that \(T\) is a \(u\)-snake of length \(1\), i.e., one \hedge of length at most two. Since \(T^{+u}\) is an \hedge of size at most one, it is a winning position for Client.
	
	Assume otherwise that \(T\) has at least two \hedges. Three cases can happen.
	\begin{itemize}
		\item \(T\) is a \(u\)-snake of length at least two. So \(T\) is a sequence of the form \((\{u,v,w\},\{v,x,y\},e_3,\ldots,e_\ell)\) where \(e_{\ell}\) has size at most \(2\) and the other \hedges have size \(3\). 
		Assume \(o(T^{+u}) = \Waiter\), then by Lemma~\ref{lem:size2hyperedge}, there is a strategy for Waiter which starts with the pair \(\{v,w\}\). However, if Client selects \(v\), then \((\{v,x,y\},e_3,\ldots,e_\ell)\) is a \(v\)-snake of size smaller than \(T\), and so the position is winning for Client by induction hypothesis. Hence \(o(T) = \Client\).
		\item \(T\) is a \(3\)-uniform simple linear \(u\)-cycle. So \(T\) is of the form \( ( \{u,v,w\},\{v,x,y\},e_3,\ldots,e_\ell) \) where \( ( \{v,x,y\},e_3,\ldots,e_\ell) \) is a \(3\)-uniform simple linear path of positive length and \(e_{\ell}\) contains \(u\).
		The conclusion is similar to the previous case. Indeed, if in \(T^{+u}\) Waiter offers the pair \(\{v,w\}\) and Client picks \(v\), we obtain a \(v\)-snake of size smaller than \(T\).
		\item \(T\) is a \(u\)-tadpole which is not a cycle. So \(T\) is a sequence of the form \((\{u,v,w\},\{v,x,y\},e_3,\ldots,e_\ell, e'_1,\ldots,e'_t)\) where \(e_{\ell}\) contains a vertex \(b\) such that \((\{u,v,w\},\{v,x,y\},e_3,\ldots,e_\ell)\) is a \(3\)-uniform simple linear \(ub\)-path, and \((e'_1,\ldots,e'_t)\) is a \(3\)-uniform simple linear \(b\)-cycle. 
		Assume \(o(T) = \Waiter\), then by Lemma~\ref{lem:size2hyperedge}, there is a strategy for Waiter which starts with the pair \(\{v,w\}\). However, if Client selects \(v\), then
		\((\{v,x,y\},e_3,\ldots,e_\ell,e'_1,\ldots,e'_t)\) is a \(v\)-tadpole of size smaller than \(T\), and so the position is winning for Client by induction hypothesis. Hence \(o(T) = \Client\). \qedhere
	\end{itemize} 
\end{proof}

We consider the set of vertices which are reachable from \(u\) by a simple linear path.
\begin{definition}
	Let \(H\) be a hypergraph and \(u\) be a vertex of \(H\). Let us denote by \(\CL_H(u)\) the induced subhypergraph of \(H^{+u}\) on the set of the vertices which are reachable from \(u\) by a simple linear path.
\end{definition}

\begin{figure}
    \centering

\begin{tikzpicture}

    \draw (0,0) node[v] (u)  {} node[above= .15cm] {$u$};
    \draw (1,-1) node[v] (a)  {};
    \draw (1,1) node[v] (b)  {};
    \draw (2,0) node[v] (v)  {} node[above= .15cm] {$v$};

    \draw (2,1.8) node[v] (c)  {};
    \draw (3,1.5) node[v] (d)  {};
    \draw (3.5,0) node[v] (e)  {};
    \draw (3,-1.5) node[v] (f)  {};
    \draw (2,-1.8) node[v] (g)  {};

\draw \hedgeiii{u}{b}{a}{5mm};
\draw \hedgeiii{b}{v}{a}{4mm};
\draw \hedgeiii{b}{c}{d}{3mm};
\draw \hedgeiii{d}{e}{f}{4mm};
\draw \hedgeiii{f}{g}{a}{3mm};

\node at (-0.2,-0.75){$e_1$};
\node at (1.5,2.1){$e_2$};
\node at (4,-0.5){$e_3$};
\node at (1.5,-2.1){$e_4$};
\node at (2.2,-0.65){$e_5$};

\end{tikzpicture}

    \caption{The vertices $u$ and $v$ are connected by the linear path $e_1e_2e_3e_4e_5$ and by the simple path $e_1e_5$ but no simple linear path connects them.}
    \label{fig:linear_path}
\end{figure}
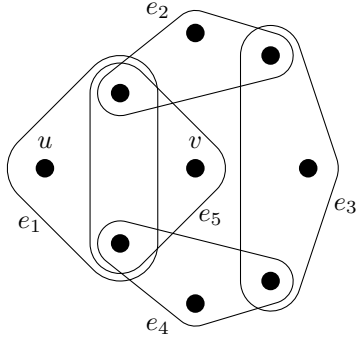

Notice that it is possible that \(v\) is reachable from \(u\) by a linear path and by a simple path with no simple linear path between \(u\) and \(v\) (see for example Figure~\ref{fig:linear_path}).

\begin{definition}
	Two vertices \(v\) and \(w\) are called {\em siblings} with respect to \(u\) if they are reachable from \(u\) by a simple linear path, but any such path to one of these vertices contain the other one. More formally, \(v \sim_u w\) if and only if \[
		v \in \CL_H(u)\setminus \CL_{H^{-w}}(u) \text{ and }w \in \CL_H(u)\setminus \CL_{H^{-v}}(u) .
	\]
\end{definition}

Notice that siblings happen only by pairs in rank \(3\) hypergraphs.
\begin{proposition}
    Let \(\hyp\) be a rank \(3\) hypergraph.
	Assume that both pairs \((v,w_1)\) and \((v,w_2)\) are siblings with respect to \(u\). Then \(w_1=w_2\).
\end{proposition}

\begin{proof}
	Let \(\mathcal{P}\) be a smallest simple linear path from \(u\) to \(v\). Since, \(v\) is not in \(\CL_{H^{-w_1}}(u) \cup \CL_{H^{-w_2}}(u)\), it means that \(w_1\) and \(w_2\) have to appear in \(\mathcal{P}\).  If \(w_1 \neq w_2\), then either \(w_1\) or \(w_2\) appears before the last \hedge in \(\mathcal{P}\). It contradicts the fact \(w_1\) and \(w_2\) do not belong to \(\CL_{H^{-v}}(u)\).
\end{proof}

\begin{definition}
	Let  \(u\) be a vertex of \(H\) hypergraph of rank \(3\). 
	We extend the hypergraph \(\CL_H(u)\) by adding, for each pair of siblings \(v \sim_u w\), a new \hedge \(\{v,w\}\). We call this hypergraph \(\CCL_H(u)\) the {\em completed} of \(\CL_H(u)\). 
\end{definition}

\begin{lemma}\label{lem:CC_iff_Tadpole}
	Let \(u\) be a vertex of \(H\) hypergraph of rank \(3\). Then, \(\uCwin\) is non empty if and only if \(o(\CCL_{H}(u)) = \Client\).
\end{lemma}

\begin{proof}
	If \(T \in \uCwin\), then \(T^{+u}\) is a subhypergraph of \(\CL_H(u) \subseteq \CCL_H(u)\). Then Lemmas~\ref{lem: sub-hypergraph} and~\ref{lem:gTadpoleCwin} imply that \(o(\CCL_H(u)) = \Client\).
	
	We prove that if \(\uCwin\) is empty then \(o(\CCL_H(u)) = \Waiter\) by induction on the number of vertices of \(\CL_{H}(u)\).  If  \(\CL_H(u)) = \emptyset\), then necessarily \(o(\CCL_H(u)) = \Waiter\).
	Otherwise, the number of vertices of \(\CL_H(u)\) is a positive integer. It means that \(\CCL_H(u)\) contains an \hedge \(e\), which is the trace of an \hedge \(\{u\}\cup e\) of \(H\). If \(e\) has size at most \(1\), then \(\{u\}\cup e\) is in \(\uCwin\). So \(e = \{v,w\}\).
	
	We have a first fact.
	\begin{fact}
		 Any vertex of \(\CL_{H^{+u-w}}(v)\) is already present in \(\CL_{H}(u)\). 
	\end{fact}
	\begin{proof}
		Let \(x\) be a vertex of \(\CL_{H^{+u-w}}(v)\). 
		It means there is an almost simple linear path \(\mathcal{P}\) in \(H^{-w}\) from \(v\) to \(x\) where only the vertex \(u\) can be repeated. 
		If \(\mathcal{P}\) contains \(u\), then, starting by its last occurrence, it contains a simple linear path from \(u\) to \(x\). Otherwise, the concatenation of \((\{u,v,w\})\) with \(\mathcal{P}\) is simple and linear. In both cases, we have that \(x\) is a vertex of \(\CL_{H}(u)\).
	\end{proof}
	
	Consequently, every vertex of \(\CCL_H(u)\) is neighbour of \(u\) or appears in some \(\CCL_{H^{+u-w}}(v)\) for such a couple of vertices \((v,w)\).
	Given such a \((v,w)\), we show, that there is a winning strategy on \(\CCL_{H^{+u-w}}(v)\) for Waiter, such that if this strategy is played on \(\CCL_H(u)^{+v-w}\), the resulting hypergraph \(G\) is a subhypergraph of \(\CCL_H(u)^{+v-w}\) and so of \(\CCL_H(u)\). The result of the lemma follows by repeating the argument for each such couple \((v,w)\) (\(G \cap \CCL_H^{+u-w'}(v')\) is a subhypergraph of \(\CCL_H^{+u-w'}(v')\) and so the strategy can be done on \(G\) by Lemma~\ref{lem: sub-hypergraph}). 
	
	So let us fix such couple \((v,w)\).
	Since \(v\) and \(w\) are vertices of \(\CL_H(u)\setminus \CL_{H^{+u-w}}(v)\), the number of vertices of \(\CL_{H^{+u-w}}(v)\) is strictly smaller than the one of \(\CL_H(u)\).

	\begin{fact}
		If \(\uCwin_{H}\) is empty, then \(\vCwin_{H^{+u-w}}\) is also empty.
	\end{fact}
	\begin{proof}
		Let us show the contrapositive. Let \(T \in \vCwin_{H^{+u-w}}\). 
		\begin{itemize}
			\item Assume first that \(T\) is reduced to one \hedge \(e\). If \(e\) is an \hedge of \(H\), then \((\{u,w,v\},e)\) is a \(u\)-snake of \(H\). Otherwise \(e\) is the trace of an \hedge \(e'\) of \(H\) of the form \(e' = \{v,x,u\}\) where \(x \neq w\) (\(e'\) can not be of length \(2\) since it would be contained in \(\{u,v,w\}\) but \(H\) is a clutter). In particular, \((\{u,w,v\},\{v,x,u\})\) is a \(3\)-uniform simple linear \(u\)-cycle of length \(2\) and so a \(u\)-tadpole in \(\uCwin\). 
			\item Assume then that \(T\) is a subhypergraph of \(H\). In particular it is a \(v\)-tadpole in \(H^{-w}\) which does not contain \(u\). Adding the \hedge \(\{u,v,w\}\) at the beginning of \(T\) gives a \(u\)-tadpole in \(\uCwin_H\).  
			\item Otherwise, \(T\) is a \(v\)-snake \((e_1,\ldots,e_{p})\) such that \(e_p\) is the trace of an \hedge \(e'= e_p \cup \{u\}\) in \(H\). 
			So adding \(\{u,v,w\}\) in front of \(T\) and replacing \(e_p\) by \(e'\) gives a \(u\)-tadpole (in fact a \(u\)-cycle) in \(\uCwin\). \qedhere
		\end{itemize}
	\end{proof}

	By the second fact \(\vCwin_{H^{+u-w}}\) is empty and so by induction hypothesis Waiter has a winning strategy \(\strat\) in \(\CCL_{H^{+u-w}}(v)\). 
	So the first fact ensures that Waiter can simulate the strategy \(\strat\) on \(H'= \CCL_{H}(u)^{+v-w} = (V',E')\). 
	After playing \(\strat\), we obtain a new hypergraph \(G\)  (depending on Client's choices). Let \(V_{\Waiter}\) and \(V_{\Client}\) be the sets of vertices claimed respectively by Waiter and Client during this phase. So \(G\) is the trace on \(V' \setminus (V_{\Waiter} \cup V_{\Client})\) of the hypergraph \(H'_{| V' \setminus V_{\Waiter}}\). We show that \(G\) is a subhypergraph of \(H'\) (i.e., there is no \hedge in \(H'_{| V' \setminus V_{\Waiter}}\) intersecting both \(V_{\Client}\) and \(V' \setminus (V_{\Waiter} \cup V_{\Client})\)).
	
	Suppose there exists \(e\) an \hedge of \(H'_{| V' \setminus V_{\Waiter}}\) such that \(x \in e \cap V_{\Client}\) and \(y \in e \setminus V_{\Client}\).
	
	\begin{itemize}
	    \item If \(e=\{x,y\}\) or if \(e=\{x,y,z\}\) with \(z\) not in \(V_{\Client}\), then any simple linear path from \(v\) to \(x\) can be extended to a simple linear path from \(v\) to \(y\) by adding the \hedge \(e\), which is impossible.
	    \item Otherwise, \(e=\{x,y,z\}\) with \(z \in V_{\Client}\). It implies that \(x\) and \(z\) are siblings with respect to \(v\) in \(H^{+u-w}\) (otherwise \(y\) would also be reachable from \(v\)). 
	    Then, there is an \hedge \(\{x,z\}\) in  \(\CCL_{H^{+u-w}}(v)\). But having \(x\) and \(z\) in \(V_{\Client}\) contradicts the fact that \(\strat\) is winning in \(\CCL_{H^{+u-w}}(v)\). \qedhere
	\end{itemize}
\end{proof}

In particular, Proposition~\ref{prop:HarderThanTadpole} is a direct consequence of Lemma~\ref{lem:CC_iff_Tadpole}.

\begin{corollary}\label{cor:CC_iff_coupleTadpoles}
	Let \(u\), \(v\) be vertices of \(H\) hypergraph of rank \(3\) such that \(\uCwin_{H^{-v}}\) and \(\vCwin_{H^{-u}}\) are empty. If \(o(H)= \Waiter\), then Waiter has a winning strategy on \(H\) which starts with the pair \(\{u,v\}\). 
\end{corollary}

\begin{proof}
	Indeed, Waiter starts offering the pair \(\{u,v\}\). Assume that Client picks \(u\) (the other case is symmetric). By hypothesis, \(\uCwin_{H^{-v}}\) is empty. By Lemma~\ref{lem:CC_iff_Tadpole} Waiter has a winning strategy \(\strat\) in \(\CCL_{H^{-v}}(u)\). Waiter can play in \(H^{+u-v}\) according to \(\strat\).
	\begin{claim}
		At the end of the strategy \(\strat\), the obtained hypergraph is a subhypergraph of \(H\).
	\end{claim}
	\begin{proofclaimitem}
		Assume that \(e\) is an \hedge of \(H^{+u-v}\) where none of its vertices are claimed by Waiter and where \(x \in e\) is a vertex claimed by Client. 
		Assume furthermore that there is \(y \in e\) which is unclaimed, i.e., it does not belong to \(\CL_{H^{-v}}(u)\).
		Two cases can happen.
		\begin{enumerate}
			\item Assume  first that \(e = \{x,y\}\) or \(e = \{x,y,z\}\) where \(z\) does not belong to \(\CL_{H^{-v}}(u)\) too. Since, \(x \in \CL_{H^{-v}}(u)\), there exists a simple linear path \(\mathcal{P}\) from \(u\) to \(x\). Adding \(e\) at the end of \(\mathcal{P}\) gives a simple linear path from \(u\) to \(y\) which is impossible.
			\item So \(x\) and \(z\) belong to \(\CL_{H^{-v}}(u)\). If one vertex of \(\{x,z\}\) is reachable from \(u\) in \(H^{-v}\) by a linear simple path which does not contain the other vertex, it would again contradicts the fact that \(y \notin \CL_{H^{-v}}(u)\). So \(x\) and \(z\) are siblings with respect to \(u\). Again, this is impossible since \(\{x,z\}\) would be a winning \hedge claimed by Client. \qedclaim
		\end{enumerate}
	\end{proofclaimitem}
	The new hypergraph is a subhypergraph of \(H\) which is winning for Waiter by hypothesis. By Lemma~\ref{lem: sub-hypergraph}, the new hypergraph is still winning for Waiter.
\end{proof}

Theorem~\ref{thm: rank3 client-waiter} follows from Corollary~\ref{cor:CC_iff_coupleTadpoles}. Indeed, it suffices to test if there is a couple of vertices \((u,v)\) such that \(\uCwin_{H^{-v}}\) and \(\vCwin_{H^{-u}}\) are empty. If it is not the case, it is a win for Client by Lemma~\ref{lem:CC_iff_Tadpole}. If such couple is found, Corollary~\ref{cor:CC_iff_coupleTadpoles} ensures, we lose nothing by starting with this couple, and we can redo the test for the smaller hypergraph.

\section{Waiter-Client games are \FPT\ on k-uniform hypergraphs} \label{sec: waiter-client}


In Oijid's thesis~\cite{oij2024thesis}, it is proved that if Waiter can win a rank~2 Waiter-Client game, she has a winning strategy in at most three moves. This leads to the following result:

\begin{proposition}[Theorem 1.70 from~\cite{oij2024thesis}]\label{prop: rank 2 Waiter-Client}
    Let $\hxf$ be a rank~$2$ hypergraph. The winner of a Waiter-Client game played on $\hyp$ can be computed in polynomial time.
\end{proposition}

We here extend this result, by proving that for any $k \ge 1$, the winner of a Waiter-Client game on a rank~$k$ hypergraph can be computed in \FPT\ time parameterized by $k$. This result is a far-reaching generalization of the following easy fact: if a rank~$k$ hypergraph has $2^k$ disjoint \hedges, it is Waiter's win.

Let $H=(V,E)$ be a $k$-uniform hypergraph. We call \emph{$\ell$-sunflower} a set $S$ of $\ell \geq 1$ \hedges of $H$ pairwise intersecting on a fixed set $C$ called \emph{center} of $S$. When the center is empty, the sunflower simply consists of disjoint \hedges. We call \emph{petal} of $S$ every set $s\setminus C$ where $s\in S$.
We authorize multisets in the definition, in particular, any \hedge $e$ of $H$ can be considered as an $\ell$-sunflower for every $\ell\geq 1$ (the emptyset being a petal). Such a sunflower with empty petals is called \emph{trivial}. The celebrated Sunflower Lemma from Erd\H{o}s and Rado~\cite{ER60} asserts that every $k$-uniform hypergraph with at least $k!(\ell-1)^k$ distinct \hedges contains a non trivial $\ell$-sunflower. We first show that the number of inclusion-wise minimal centers is bounded in terms of $k$ and $\ell$. This is the first step of the FPT algorithm for the Waiter-Client game.

We say that an $\ell$-sunflower $S$ of $H$ is \emph{minimal} (with respect to inclusion) if 
no $\ell$-sunflower of $H$ has a center strictly included in the center of $S$. Let $Y$ be a subset of vertices of $H$, we say that $S$ is \emph{outside} $Y$ if all petals of $S$ are disjoint from $Y$. In particular, every \hedge \(e \in E\) forms a trivial \(\ell\)-sunflower outside \(Y\) for every subset \(Y \subseteq V\).

\begin{lemma}\label{lem:minsun}
There exists a function $mc_k$ for which every $k$-uniform hypergraph $H$ has at most $mc_k(\ell)$ distinct centers of minimal $\ell$-sunflower. Moreover, there exists a function $omc_k$ such that whenever $Y$ is a subset of vertices of size $y$, $H$ has at most $omc_k(\ell,y)$ distinct centers of minimal $\ell$-sunflower outside $Y$.
\end{lemma}

\begin{proof}
We focus on the existence of $mc_k$, and postpone the outside case $omc_k$ to the end of the proof. 

We just have to show that if $C_1,\dots ,C_t$ are distinct centers of size $c$ of a family of minimal $\ell$-sunflowers $S_1,\dots ,S_t$, then $t$ is bounded in terms of $k$ and $\ell$. Indeed, this implies the first part of the lemma since \(mc_k(\ell)\) will be at most \(k\) times this bound. Let us fix $t'=2k(\ell-1)+1$.
We denote by $P_i^j$ the $j^{th}$ petal of $S_i$. If  $t$ is large enough, we can extract a subfamily of $t'$ sunflowers from $S_1,\dots ,S_t$ which, free to reorder, is assumed to be $S_1,\dots ,S_{t'}$ with the following additional properties:

\begin{itemize}
    \item all (distinct) centers $C_1,\dots ,C_{t'}$ form a sunflower with center $X$. 
    \item for all $1\leq j\leq \ell$, all petals $P_1^j,\dots ,P_{t'}^j$ form a sunflower with center $Q_j$. 
\end{itemize}

Indeed, we just have for this to iterate $\ell +1$ extractions using the Sunflower Lemma, one for the centers, and $\ell$ for the petals. To conclude, we now extract an $\ell$-sunflower $S$ centered at $X$ in \(H\). This is indeed a contradiction to the minimality of centers $C_i$. A simple greedy argument suffices for this.

Let us say that $B_i:= C_i\setminus X$ and $D_i^j:=P_i^j\setminus Q_j$. With this notation, the $j^{th}$ \hedge of the sunflower $S_i$ is $X\cup B_i\cup Q_j\cup D_i^j$. 

Our first \hedge $e_1$ of $S$ is the first \hedge of $S_1$, that is $C_1\cup P_1^1$. Observe that $e_1$ intersects at most $k$ of the (disjoint) subsets $B_i$ and also at most $k$ of the subsets $D_i^2$. Since $t'>2k$, there is an index $a$ such that $e_1$ is disjoint from $B_a\cup D_a^2$ (notice that \(a\) is automatically distinct from \(1\) since \(e_1\) is disjoint of \(B_a\)). We pick the \hedge $e_2$ as the second \hedge of $S_a$, that is $X\cup B_a\cup Q_2\cup D_a^2$. Since all \hedges of $S_1$ intersect on $C_1$, the \hedge $e_1$ is disjoint from $Q_2$. In particular $e_1\cap e_2=X$. By a similar argument, since $t'>4k$, there is an index $b$ such that $e_1\cup e_2$ is disjoint from $B_b\cup D_b^3$. We now pick the \hedge $e_3$ as the third \hedge of $S_b$. Iterating the process leads to an $\ell$-sunflower $S$ centered on $X$, a contradiction.

The proof for sunflowers outside some set $Y$ is similar, but we have to take care of the fact that the final sunflower $S$ must have its petals outside $Y$. This is granted for the subsets $Q_j$ and $D_i^j$, which are petals of the original $\ell$-sunflowers $S_1,\dots ,S_t$ outside $Y$. But the (disjoint) sets $B_i$ can intersect $Y$. Observe that this happens at most $|Y|$ times. Hence, choosing $t'=|Y|+2k(\ell-1)+1$ leads to a contradiction.
\end{proof}


We now turn to the key-definition. Given some integer $\ell$, we say that a set $K\subseteq V$ is an \emph{$\ell$-kernel} of $H$ if for every \hedge $e\in E$, there exists $C\subseteq e\cap K$ and an $\ell$-sunflower $S$ outside $K$ centered at $C$. Note that we can assume that $S$ is minimal outside $K$. 
Observe that the union of all \hedges of $H$ is an $\ell$-kernel for every $\ell$, and that if $H$ has $\ell$ disjoint \hedges, then $\emptyset$ is an $\ell$-kernel. We now show that there is always a bounded size kernel.

\begin{theorem}\label{th:kernel}
There exists a function $f_k$ such that every $k$-uniform hypergraph $H$ has an $\ell$-kernel of size at most $f_k(\ell)$, for every integer $\ell$.
\end{theorem}

\begin{proof}
Let $X_0=\emptyset$ and $\ell$ be some integer. We denote by $C_0$ the set of all centers of minimal $\ell$-sunflowers, which is bounded by $mc_k(\ell)$ by Lemma~\ref{lem:minsun}. We then set $X_1:=\bigcup C_0$ and start to define our sequence $(X_i)$ as follows: Once $X_i$ is defined, we denote by $C_i$ the set of all centers of minimal $\ell$-sunflowers outside $X_i$ and set 
$X_{i+1}:=(\bigcup C_i)\cup X_i$. Our goal is to show that there exists a step $i$, bounded by a function on $k$ and $\ell$, such that $X_{i+1}=X_i$. Observe that when this step is reached, $X_i$ is an $\ell$-kernel. Indeed, let $e\in E$ be an \hedge. Note that $e$ by itself is a trivial $\ell$-sunflower outside $X_i$, hence there is a $C$ included in $e$ which is the center of a minimal $\ell$-sunflower outside $X_i$. By definition of $X_{i+1}$, we have $C\subseteq X_{i+1}=X_i$, and therefore $C\subseteq X_i\cap e$, implying that $X_i$ is an \(\ell\)-kernel.

To reach our conclusion, we need to show that the size of $X_i$ is upper bounded by some fixed function $g(i)$ (depending on $k$ and $\ell$). For this, we just have to bound the size of $X_{i+1}$ in terms of the size of $X_i$. This is easily obtained by Lemma~\ref{lem:minsun}, since the total number of minimal $\ell$-sunflowers outside $X_i$ is at most $omc_k(\ell,|X_i|)$. 

 Let us now denote, for every $i$, the sequence $s_i=(t_0(i),t_1(i),\dots ,t_{k}(i))$ where $t_c(i)$ denotes the number of distinct centers $C$ of size $c$ of minimal $\ell$-sunflowers outside $X_i$. If $X_i\neq X_{i+1}$, we claim that $s_{i+1}$ is lexicographically smaller than $s_i$. To see this, assume that $c$ is minimum such that $t_c(i)\neq t_c(i+1)$ and suppose for contradiction that $t_c(i)<t_c(i+1)$.
 %
 %
 %
 This means that there is a set $C$ of size $c$ which is the center of a minimal $\ell$-sunflower outside $X_{i+1}$, but not the center of a minimal $\ell$-sunflower outside $X_{i}$. In particular, $C$ strictly contains $C'$ which is the center of a minimal $\ell$-sunflower outside $X_{i}$, but not the center of a minimal $\ell$-sunflower outside $X_{i+1}$ (by minimality of $C$). Since $t_{|C'|}(i)=t_{|C'|}(i+1)$, there is a set $C''$ of size $|C'|$ which is the center of a minimal $\ell$-sunflower outside $X_{i+1}$, but not the center of a minimal $\ell$-sunflower outside $X_{i}$. This process always find a smaller center which is minimum outside $X_{i+1}$ but not $X_i$. This is a contradiction.

 Note that the previous argument moreover shows that if $s_i=s_{i+1}$, then the centers of minimal $\ell$-sunflowers outside $X_i$ and outside $X_{i+1}$ are the same, thus $X_i=X_{i+1}$. 
 
 Therefore, we just have to show that the sequence $(s_i)$ is ultimately constant after some fixed number of steps. 
 Let us examine for this how the sequence $(s_i)$ evolves from $i$ to $i+1$. At each step, the minimum index coordinate $t_c$ which differs in $s_i$ and $s_{i+1}$ decreases. Meanwhile, the coordinates $t_{c'}$ where $c'>c$ can change in an arbitrary way, but no more than some fixed amount since the size of $X_{i+1}$ is bounded by $g(i+1)$ (a crude upper bound for the increase of $t_{c'}$ would be for instance $omc_k(\ell,g(i+1))$). In all, this process terminates before some fixed number of steps depending on $k$ and \(\ell\). 
\end{proof}

Kernels are relevant for Waiter-Client games. Indeed, if $K$ is a $2^k$-kernel of a $k$-uniform hypergraph $H=(V,E)$ and $T_K(H)$ is the trace hypergraph on vertex set $K$ and \hedge set $E_K=\{e\cap K:e\in E\}$, we have the following equivalence:

\begin{lemma}\label{lem:equiv}
$H$ is Waiter win if and only if $T_K(H)$ is Waiter win. 
\end{lemma}

\begin{proof}
Since $T_K(H)$ is obtained by reducing the size of some \hedges of $H$ (and deleting isolated vertices), if Waiter has a winning strategy on $H$, it has one on $T_K(H)$.

Now assume that $T_K(H)$ is Waiter win. Waiter can play her strategy on $T_K(H)$ in order to select a set of vertices which contains some $e\cap K$ for $e\in E$. Since $K$ is a kernel, there exists a $2^k$-sunflower $S$ outside $K$ with center $C\subseteq e\cap K$. Now Waiter just has to play outside $K$ to select one of the $2^k$ petals of $S$ and win the game.
\end{proof}

We can now prove Theorem~\ref{th:FPT}, asserting that Waiter-Client is FPT on $k$-uniform hypergraphs. 

\begin{proof}
Let $H$ be some input $k$-uniform hypergraph on $n$ vertices and $m$ \hedges. We first describe an algorithm running in time $O(h(k).(n+m)^2)$ which decides if $H$ is Waiter-win. 

The strategy is to compute a $2^k$-kernel $K$ as in Theorem~\ref{th:kernel} which size only depends on $k$. Then, by Lemma~\ref{lem:equiv}, we just have to (brute force) test in $O(BF(k))$ time if $H_K$ is a Waiter win. The only point to check is how efficiently one can compute the sequence $(X_i)$ as in the proof of Theorem~\ref{th:kernel}.

As the argument is similar for getting $X_{i+1}$ from $X_i$, we just describe how to compute $X_1$, that is the set of centers of minimal $2^k$-flowers. An easy algorithm consists in (bottom up) testing  all $2^k m$ subsets included in some \hedge of $H$ and determine which ones are centers of $2^k$-sunflowers (this single test being done in linear time, the computation of all centers can be achieved in quadratic time). Once computed, determining which centers are minimal can be done in linear time by dynamic programming. In all, each $X_i$ can be computed in quadratic time, and the number of steps is bounded in terms of $k$.

To improve the complexity of the above algorithm to linear, we need to show that, while reading all \hedges of $H$, we can only keep constant size information for each possible center $C$ to determine if it is the center of a $2^k$-sunflower. For this, observe that whenever we have read $k!(k2^k-1)^k$ \hedges containing $C$, then there is a non trivial $k2^k$-sunflower $S$ centered at $C'$ containing $C$. Now we can replace all \hedges in $S$ by $C'$ since the existence of a $2^k$-sunflower centered at $C$ is left unchanged. We can then compute in linear time all the centers $C$, and proceed as previously.
\end{proof}


\section{Open problems} \label{sec: open}

A reasonable guess is that Waiter-Client game is NP-hard when the set of \hedges is explicit, i.e. the size of the input is the sum of the sizes of \hedges. This is our first open problem. This would highlight the fact that Fixed Parameter Tractability is probably the best one can expect, leaving open the order of magnitude of the complexity in $k$. Observe that the proof of Theorem~\ref{th:FPT} shows in particular that the minimum number of moves in an optimal winning strategy for Waiter is at most some value $os(k)$. Proposition~\ref{prop: rank 2 Waiter-Client} gives $os(2)=3$. Alas, we have no decent bound to propose for $os(3)$, and a very optimistic analysis of the bound provided by Theorem~\ref{th:kernel} already gives an Ackermann type bound for $os(k)$. We believe that a much more reasonable upper bound can be achieved for $os(k)$.

The relationship between positional games and hypergraph colorability has driven a lot of research in Maker-Breaker. This is also the case here: 
Is it possible that Waiter's win in Clien-Waiter implies that the hypergraph is \(2\)-colorable?

\bibliography{bibliography} 
\bibliographystyle{alpha} 

\end{document}